\documentclass{amsart}
\usepackage{amsmath, amssymb, amsthm, graphics}
\theoremstyle{plain}
\newtheorem{prop}{Proposition}[section]
\newtheorem{coro}[prop]{Corollary}
\newtheorem{lemm}[prop]{Lemma}
\newtheorem{algo}[prop]{Algorithm}

\theoremstyle{definition}
\newtheorem*{conv}{Convention}
\newtheorem{defi}[prop]{Definition}

\newtheorem{exam}[prop]{Example}
\newtheorem{ques}[prop]{Question}
\newtheorem{rema}[prop]{Remark}

\numberwithin{equation}{section}
\numberwithin{figure}{section}

\def\Refff#1; #2; #3; #4; #5; #6; #7\par{%
\bibitem{#1} #2, {\it #3}, #4 {\bf #5} (#6) #7}

\def\Reff#1; #2; #3; #4\par{%
\bibitem{#1} #2, {\it #3}, #4}

\def\Ref#1; #2; #3; #4\par{%
\bibitem{#1} #2, #3, #4}
\def\;{\, ; \,}

\def\cl#1{\overline{\vrule width0pt height5pt #1}}

\let\D=\Delta
\def\divel{\mathrel{\preceq_l}}
\def\divl{\mathrel{\prec_l}}
\def\divr{\mathrel{\prec_r}}
\def\dr{\mathord{\backslash}}

\let\e=\varepsilon
\def\eqp{\equiv}
\def\eqpm{\equiv^{\scriptscriptstyle \!\pm}}

\def\GR(#1; #2){\langle #1 \; #2 \rangle}

\def\ie{{\it i.e.}}
\def\ii {^{-1}}
\def\ince{\subseteq}

\def\lg(#1){{\rm lg}(#1)}
\let\ll=\lambda

\def\MS{M_{S}}
\def\MO(#1; #2){\langle #1 \; #2
\rangle^{\scriptscriptstyle\!+}}
\def\multl{\succ_l}
\def\multr{\succ_r}

\let\o=\omega
\def\op{\dot+}

\let\pp=\dots

\let\r=\rho
\def\resp{\hbox{resp{.} }}
\def\RR{{\mathcal R}}		
\def\rvl{\curvearrowright_l}
\def\rvlq{\curvearrowright_l^{\scriptscriptstyle\sharp}}
\def\rvr{\curvearrowright_r}
\def\rvrq{\curvearrowright_r^{\scriptscriptstyle\sharp}}
\def\rv{\curvearrowright}

\let\s=\sigma
\def\SS{{\mathcal S}}
\def\SSp{\SS'}
\def\SSr{\widehat\SS}

\def\xx{x}		
\def\yy{y}
\def\zz{z}

\def\ss{s} 
\def\tt{t}
\def\rr{r}

\def\uu{{\mathbf u}}
\def\vv{{\mathbf v}}
\def\ww{{\mathbf w}}

\begin{document}

\author{Patrick DEHORNOY}
\address{Laboratoire SDAD, Math\'ematiques\\
Universit\'e de Caen BP 5186, 14032 Caen, France}
\email{dehornoy@math.unicaen.fr}
\urladdr{//www.math.unicaen.fr/\!\!\!\hbox{$\sim$}dehornoy}

\title{COMPLETE POSITIVE GROUP PRESENTATIONS}



\begin{abstract}
 A combinatorial property of prositive group
 presentations, called completeness, is introduced,
 with an effective criterion for recognizing complete
 presentations, and an iterative method for
 completing an incomplete presentation. We show
 how to directly read several properties of the
 associated monoid and group from a complete
 presentation: cancellativity or existence of common
 multiples in the case of the monoid, or isoperimetric
 inequality in the case of the group. In particular, we
 obtain a new criterion for recognizing that a monoid
 embeds in a group of fractions. Typical presentations
 eligible for the current approach are the standard
 presentations of the Artin groups and the Heisenberg
 group.
\end{abstract}

\maketitle

\section*{Introduction}

This paper is about monoids and groups defined by a
presentation. As is well-known, it is hopeless to directly
read from a presentation the properties of a group or a
monoid: even recognizing whether the group is trivial
is undecidable in general~\cite{Mil}. However,
partial results may exist when one restricts to
presentations of a special form: a typical example is the
small cancellation theory, in which a number of
properties are established for those groups or monoids
defined by presentations satisfying some conditions
about subword overlapping in the relations \cite{Hig,
Kas, LyS, Sch}. Another example is Adyan's criterion
\cite{Ady, Rem} which shows that a presented monoid
embeds in the corresponding group if there is
no cycle in some graph associated with the
presentation. The aim of this paper is to study a
combinatorial property of positive group
presentations (\ie, of presentations where all
relations are of the form $u = v$ with only positive
exponents in~$u$ and~$v$) that we call completeness,
and to show that several nontrivial properties of the
associated monoid and group can be read directly
when a complete presentation is known: the
properties we shall investigate here are cancellativity,
existence of common multiples, embeddability in a
group of fractions in the case of the monoid, solution
for the word problem, and isoperimetric inequality in
the case of the group. What we do in each case is to
give sufficient conditions for the monoid or the group
defined by a supposedly complete presentation to
satisfy the considered property. A typical example is
Prop.~\ref{P:canc}, which states that, if $(\SS, \RR)$ is
a complete presentation, then a sufficient condition for
the associated monoid to be cancellative is that
$\RR$ contains no relation of the form
$\ss u = \ss v$  or $u \ss = v \ss$ with $u \not= v$:
thus, if there is no obvious counter-example to
cancellativity, then there is no hidden
counter-example either.

The interest of such results could be void if complete
presentations did not exist. Actually, they do: it is
even trivial that every group admits complete
presentations---as the name suggests, a complete
presentation is one with enough relations, and the full
presentation consisting of all relations is always
complete. The interesting case is when
there exists a finite (or, at least, simple) complete
presentation: we shall see that this happens for a
number of groups, such as many generalized braid
groups (in particular some of those associated with
complex reflection groups~\cite{BMR}), more generally
all Garside groups of~\cite{Dgk}, but also quite
different groups, such as the Heisenberg group, which
is nilpotent.

The main technical ingredient we shall use is a
combinatorial transformation called {\it word
reversing}. It is a refinement of the monoid
congruence, in the sense that applying reversing to a
word gives an equivalent word, but, in general, the
converse is not true, \ie, it is not true that any pair of
equivalent words can be produced (or, better,
detected) using reversing. Essentially, we say that a
presentation is complete when the latter occurs, in
which case the uneasy study of word equivalence can
be replaced with the easier study of reversing. 

It seems that the reversing process has been first
considered in~\cite{Dfa}, and it has been
investigated---and in particular some notion of
completeness has been considered---in several
papers~\cite{Dff, Dfx, Dgc, Dgd}, but so far always in
the particular case of presentations with few relations,
namely the so-called complemented presentations
where there exists at most one relation
$\ss \cdots = \tt \cdots$ for each pair of letters~$\ss$,
$\tt$. K.~Tatsuoka in~\cite{Tat} (in the case of Artin
groups) and R.~Corran in~\cite{Cor} (in the case of
singular Artin monoids) have independently developed
equivalent processes in slightly different frameworks,
but always with equally or more restricted initial
assumptions. 

The current work addresses arbitrary positive
presentations. The advantage of such a
generalization---which forces to renew the technical
framework---does not only lie in the new
groups that become eligible, but
rather in the underlying change of viewpoint.
Previously, the principle was to study the possible
completeness of a (complemented) presentation: in
good cases, the presentation was complete and one
could deduce consequences---as in the case of the
standard presentation of the braid groups \cite{Gar} or
of their alternative presentation
of~\cite{BKL}---otherwise, if the presentation was not
complete, one could say nothing. Our current
approach enables us not only to establish the
completeness of a presentation, but also, if needed, to
complete an initially incomplete presentation. This
completion process may require an infinite number of
steps, but, in good cases, it is a finite procedure, and we
shall see on examples how it enables us to investigate
some monoids or groups that remained outside the
range of all previously known methods. In particular,
we obtain a new method for proving that a monoid
embeds in a group of fractions, and apply it to answer a
question of~\cite{HaK} about a nonstandard
presentation of Artin's braid group~$B_3$ introduced
by V.~Sergiescu in~\cite{Ser}. 

One of the applications of word reversing is (in
good cases) a solution of the word problem. Let us
mention here some similarity between this solution and
Dehn's algorithm for hyperbolic groups: in both cases,
the idea is to decide whether a word
represents~$1$ without introducing any new pair of
generators~$\ss\ss\ii$ or~$\ss\ii\ss$. However,
contrary to Dehn's algorithm, the reversing algorithm
may increase the length of the words, and it is not
linear in general, but, on the other hand, it works for
groups that are not word hyperbolic, such as the braid
groups, or even the Heisenberg group, whose
isoperimetric function is known to be cubic.

The rather vague description above might remind the
reader of the Knuth-Bendix completion
method~\cite{KnB, Coh}, which also consists in
starting with a group presentation, possibly adding
some consequences of the initial relations, and
obtaining a so-called complete rewrite-system that
enables one to solve the word problem---see
\cite{HeM} for examples in the case of spherical Artin
groups. The similarity with the current approach is
superficial only: our method also possibly provides a
solution to the word problem by means of rewriting
rules, and the r\^ole of the cube conditions in our
completion procedure is analogous to that of critical
pairs in~\cite{KnB}, but there seems to be no more
precise connection in general, and we do not see how to
attach any confluent rewrite-system to the
combinatorial word transformations we consider, in
particular because we simultaneously use positive
and arbitrary words, \ie, we work both with the
monoid and the group. Actually, more than in the
Knuth-Bendix method, our approach originates in
Garside's analysis of the braid monoids~\cite{Gar}:
with our current definitions, the proof of
Prop.~H in~\cite{Gar}, as well as that of the
K\"urzungslemma of~\cite{BrS} is a proof that the
standard presentations of the (generalized) braid
groups is complete.

The paper is organized as follows. In Sec.~\ref{S:redr},
we define the general reversing process and establish
its basic properties. Then, in Sec.~\ref{S:comp}, we
introduce completeness, and, again, establish basic
results, in particular that every monoid admits a
complete presentation. In Sec.~\ref{S:cube}, we
introduce the cube condition, a technical property
which we show is equivalent to completeness. We use it
to establish our main criterion for recognizing
completeness in~Sec.~\ref{S:crit} and, in
Sec.~\ref{S:completion}, to complete initially
incomplete presentations. The rest of the paper is
devoted to studying monoids and groups from a
complete presentation. In Sec.~\ref{S:mono}, we
consider properties of the monoid: cancellativity, word
problem, common multiples. Finally, in
Sec.~\ref{S:embe}, we investigate similar questions for
the group: recognizing groups of fractions, solving the
word problem, computing bounds for the isoperimetric
function.

\begin{conv}
 A number of notions will appear with a right and
 a left version. We shall use $r$- for ``right'' and $l$- for
 ``left'': $r$-reversing, $r$-completeness, etc.
\end{conv}

\section{Reversing}\label{S:redr}

Our aim is to study groups and monoids from a
presentation. Here we consider positive group 
presentations, defined as those presentations where all
relations have the form $u = v$, where $u$ and $v$
are nonempty positive words, \ie, inverses of the
chosen generators do not occur in~$u$ or~$v$. At the
expense of adding new generators, this is not a
restriction in the case of groups, but this means that we
restrict to monoids with non nontrivial units. Our
notation will be as follows. If $\SS$ is a nonempty set,
we denote by~$\SS^*$ the free monoid generated
by~$\SS$, \ie, the set of all words on~$\SS$ equipped
with concatenation; we use $\e$ for the empty word. A
positive group presentation is then a pair~$(\SS, \RR)$
where $\RR$ is a family of pairs of nonempty words
in~$\SS^*$, the relations of the presentation. As usual,
we shall often write $u = v$ instead of $\{u, v\}$ for a
relation. We denote by~$\MO(\SS; \RR)$ the monoid
associated with the presentation~$(\SS, \RR)$,
\ie, the monoid~$\SS^* \! / \! \eqp$, where $\eqp$
is the smallest congruence on~$\SS^*$ that
includes~$\RR$. Then, we denote by~$\GR(\SS; \RR)$
the associated group: introducing for each letter~$\ss$
in~$\SS$ a disjoint copy~$\ss\ii $ of~$\ss$, and
using~$\SS\ii $ for the set of all~$\ss\ii $'s, the group
$\GR(\SS; \RR)$ is $(\SS \cup \SS\ii )^* \! / \! \eqpm$, where $\eqpm$ is
the smallest congruence on $(\SS \cup
\SS\ii )^*$ that includes~$\RR$ (hence~$\eqp$) and
contains all pairs $\{\ss \ss\ii , \e\}$, $\{x\ii x, \e\}$,
\ie, all relations $\ss\ss\ii =
\ss\ii \ss = \e$, for~$\ss$ in~$\SS$. For
$\ww$ a word on~$\SS \cup \SS\ii $, we denote
by~$\ww\ii $ the word obtained from~$\ww$ by
exchanging~$\ss$ and~$\ss\ii $ everywhere and
reversing the order of the letters: if $\ww$
represents~$\xx$ in~$\GR(\SS; \RR)$, then $\ww\ii $
represents~$\xx\ii $.

\begin{conv}
 In the previous framework, we reserve $\ss$, $\tt$
 for letters in~$\SS$, and $u$, $v$, $w$ for words
 in~$\SS^*$. We use bold letters $\uu$,
 $\vv$, $\ww$ for words on the symmetrized
 alphabet~$\SS \cup \SS\ii $. For $u$ a word
 in~$\SS^*$, we shall use $\cl u$ for the element of
 the considered monoid~$\MO(\SS; \RR)$ represented
by~$u$.
\end{conv}

Our main tool in the sequel is a combinatorial
transformation of words called reversing.

\begin{defi}\label{D:redr}
 Assume that $(\SS, \RR)$ is a positive group
 presentation, and $\ww$, $\ww'$ are words on~$\SS
 \cup \SS\ii $. We say that  $\ww \rvr^{(1)} \ww'$ is
  true if $\ww'$ is obtained from~$\ww$ 
 
 - either by deleting some subword $u\ii u$ where
 $u$ is a nonempty word on~$\SS$, 
 
 - or by replacing some subword $u\ii v$ where $u$,
 $v$ are nonempty words on~$\SS$ with a word~$v'
 {u'}\ii$ such that $u v' = v u'$ is a relation
 of~$\RR$. 
 
 Defining an {\it $r$-reversing
 sequence} to be a (finite or infinite) sequence of
 words $\ww_0, \ww_1, \pp$ satisfying $\ww_i
 \rvr^{(1)} \ww_{i+1}$ for every~$i$, we write
 $\ww \rvr^{(k)} \ww'$ if there exists a
 length~$k$ $r$-reversing sequence
 from~$\ww$ to~$\ww'$, and we say that $\ww$ is
 {\it $r$-reversible} (\ie, right reversible)
 to~$\ww'$---or that $\ww$ reverses to~$\ww'$ on the
 right---denoted $\ww \rvr \ww'$ if $\ww \rvr^{(k)}
 \ww'$ holds for some nonnegative integer~$k$.
 
 Symmetrically, we say that $\ww$ is {\it
$l$-reversible}
 to~$\ww'$, denoted $\ww \rvl \ww'$, if $\ww'$ is obtained
 from~$\ww$ by repeatedly deleting
 subwords~$u u\ii $ and
 replacing subwords~$u v\ii $ with
 words~${v'}\ii u'$ such that $v' u = u' v$ is a
 relation of~$\RR$. 
\end{defi}

\begin{figure}[htb]
 $$\includegraphics{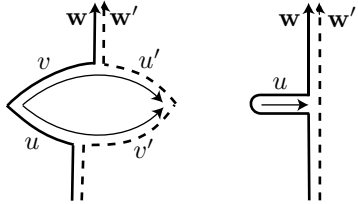}$$
 \caption{Right reversing in the Cayley graph}
 \label{F:rvca}
\end{figure}

Fig.~\ref{F:rvca} illustrates reversing in the Cayley
graph of~$\GR(\SS; \RR)$: a relation $u v' = v u'$
corresponds to an oriented cell, and the words~$\ww$,
$\ww'$ correspond to paths; then saying that $\ww
\rvr^{(1)} \ww'$ is true means that the path associated
with~$\ww'$ is obtained from that associated
with~$\ww$ by reversing the way the cell $u v' = v u'$
is crossed, namely going through the final vertex
instead of through the initial one. The case when we
delete $u\ii u$ is not particular provided we assume
that the trivial relation $u = u$ is added to the
presentation.

The study of $l$-reversing is of course similar to that
of $r$-reversing. However the reader should keep in
mind that $u v\ii \rvr v' {u'}\ii $ does {\it not} imply
$v' {u'}\ii \rvl u v\ii $: deleting $\ss\ii 
\ss$ is not a reversible process, and we always have $\ss\ii 
\ss \rvr \e$, but never $\e \rvl \ss\ii \ss$. 

\begin{exam}\label{X:aabb}
 Consider the presentation 
 $(a, b ; a^2 = b^2, ab = ba)$,
 and let $\ww = a\ii b a b\ii $. By
 using the first relation, we find $\ww \rvr^{(1)} b a\ii 
 a b\ii \rvr^{(1)} bb\ii $, hence $\ww \rvr^{(2)} b
 b\ii $, and no further $r$-reversing is possible. By
 using the second relation
 first, we can construct a different $r$-reversing
 sequence, for instance $\ww \rvr^{(1)} a b\ii a
 b\ii \rvr^{(1)} a^2 b^{-2}$. Observe that the
 previous sequences are maximal in the sense that
 they end up with a word of the form~$v u\ii $ with
 $u$, $v$ in~$\SS^*$, and no further $r$-reversing is
 possible as such a
 word contains no subword of the form~${u'}\ii v'$
 with $u'$, $v' \not= \e$.
 An example of a (maximal) $l$-reversing sequence
 is $\ww \rvl^{(1)} a\ii b b\ii a
 \rvl^{(1)} a\ii a$.
\end{exam}

As the previous example shows, reversing is not a
deterministic process in general: there can exist many 
ways of reversing one word. The only case where
$r$-reversing is certainly deterministic is the case of
complemeneted presentations:

\begin{defi}
  A positive presentation~$(\SS, \RR)$ is said to be
  {\it $r$-complemented\/} if, for all letters~$\ss$,
  $\tt$ in~$\SS$, there is at most one relation of the
  type $\ss \cdots = \tt \cdots$ in~$\RR$, and no
  relation of the type $\ss \cdots = \ss\cdots$. We say
  that $(\SS, \RR)$ is {\it complemented\/} if it is both
  $r$- and $l$-complemented, the latter being defined
  symmetrically.
\end{defi}

Reversing has been investigated in the complemented
case in~\cite{Dff} and~\cite{Dgc}. The purpose of our
current study is to extend the results to the general
case, \ie, to non necessarily complemented
presentations. We hope to convince the reader that
this extension is not trivial and that the general case is
actually the most convenient one, in particular because
it forces us to carefully choose the right technical
conditions whereas an additional superfluous
hypotheses like complementedness left some
misleading flexibility.

\begin{figure}[htb]
 \includegraphics{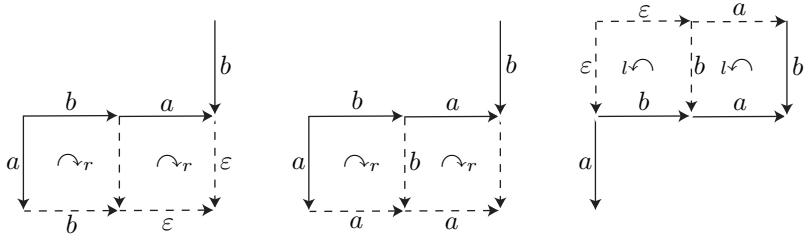}
 \caption{Two $r$- and one $l$-reversing
 sequences from~$a\ii bab\ii $}\label{F:exre}
\end{figure}

It is convenient to associate with every $r$-reversing
sequence $\ww_0, \ww_1, \pp$ a labelled planar
graph as follows. First, we associate with~$\ww_0$ a
path labelled with the successive letters
of~$\ww_0$: we associate to every positive letter~$\ss$
an horizontal right-oriented edge labelled~$\ss$, and
to every negative letter~$\ss\ii $ a vertical
down-oriented edge labelled~$\ss$. Then we by and by
represent the words~$\ww_1$, $\ww_2$,
\pp as follows: if $\ww_{i+1}$ is obtained from~$\ww_i$ by
replacing~$u\ii v$ with~$v' {u'}\ii $ (such that $u v' =
v u'$ is a relation of our presentation), then the
involved factor~$u\ii v$ is associated with a
diverging pair of edges in a path labelled~$\ww_i$
and we complete our graph by closing the open
pattern~$u\ii v$ using horizontal edges labelled~$v'$
and vertical edges labelled~$u'$:
\[
 \includegraphics{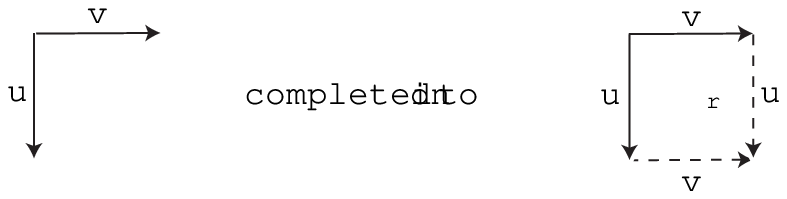}
\]
The case of the empty word~$\e$, which appears
when a factor $u\ii u$ is deleted or some relation
$u v' = v$ is used, is treated similarly: we introduce
$\e$-labelled edges and use them according to the
conventions $\e\ii u \rvr u \e\ii $, $u\ii \e \rvr \e
u\ii $, and $\e\ii \e \rvr \e \e\ii $. A symmetric
construction is associated with $l$-reversing. With
these conventions, the graphs associated with the
reversing sequences of Example~\ref{X:aabb} are those
represented in Fig.~\ref{F:exre}. Notice that the
reversing graphs, which are reminiscent of van
Kampen diagrams, need not be fragments of the
Cayley graph: several vertices may represent the same
element of the group, and they are not identified.

Let us turn to the technical study of reversing. First,
we observe that we can restrict without loss of
generality to reversing transformations of a particular
type, namely those involving length~$2$ initial
factors, \ie, to the case when $u$ and $v$ are single
letters.

\begin{lemm}\label{L:rest}
 Let $\rvr'$ be the binary relation defined as~$\rvr$
 excepted that we require that the words~$u$
 and~$v$ have length~$1$ exactly. Then $\rvr'$
 coincides with~$\rvr$.
\end{lemm}

\begin{figure}[htb]
 $$\includegraphics{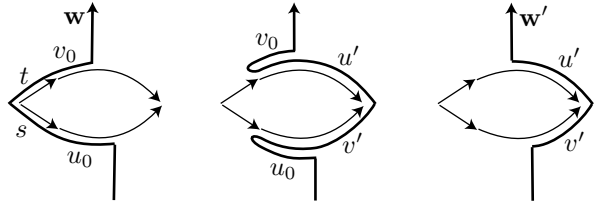}$$
 \caption{Restricted reversing}
 \label{F:rest}
\end{figure}

\begin{proof}
 (Fig.~\ref{F:rest}) By definition, $\rvr'$ is included
 in~$\rvr$. So it suffices that we prove that $\ww
 \rvr^{(1)} \ww'$ implies $\ww \rvr' \ww'$. Assume that
 $u v' = v u'$ is a relation of~$\RR$, with $u$, $v \not=
 \e$. Let $\ss$ and $\tt$ be the first letters
 of~$u$ and~$v$, say $u = \ss u_0$ and $v = \tt
 v_0$. By hypothesis, $\ss u_0 v' = \tt v_0 u'$ is a
 relation of~$\RR$, so we find
 \[
 u\ii v = u_0\ii \ss\ii \tt v_0
 \rvr' u_0\ii u_0 v' {u'}\ii v_0\ii v_0
 \rvr' v' {u'}\ii,
 \]
 as, by construction, $w\ii w \rvr' \e$ holds for every
 word~$w$ in~$\SS^*$.
\end{proof}

\begin{rema}\label{R:word}
 Instead of restricting the definition of reversing by
 considering particular subwords~$u\ii v$, we can
 extend it by relaxing the assumption that $u$
 and~$v$ are nonempty. Merely dropping
 the assumption would allow one to replace~$\e$
 by any word~$u v\ii$ such that $u = v$ is a relation
 of~$\RR$, which contradicts the implicit underlying
 principle that reversing should not increase
 complexity. But an interesting notion is obtained when
 we allow $u$ to be empty provided $u'$ is empty as
 well, \ie, we allow replacing~$v$ with~$v'$ when $v =
 v'$ is a relation of~$\RR$, and, symmetrically, we
 allow  $v$ to be empty provided $v'$ is, \ie, 
 we allow replacing $u\ii$ with~${u'}\ii$ when $u = u'$
 is a relation of~$\RR$. Most of the subsequent study
 of~$\rvr$ remains valid when the extended
 relation~$\rvrq$ so defined replaces~$\rvr$.
 However, in pratice, in particular when
 implementations are concerned, using~$\rvrq$
 instead of~$\rvr$ makes the verifications longer, as
 more transformations have to be considered. 
\end{rema}

We establish now some general properties of
(right) reversing. Owing to Lem\-ma~\ref{L:rest}, we
can always assume without loss of generality that the
basic reversing steps involve factors of the form~$\ss\ii
\tt$ where $\ss$ and~$\tt$ are single letters.

\begin{lemm}\label{L:symm}
 For all words~$\ww$, $\ww'$ on~$\SS \cup \SS\ii$,
 $\ww \rvr \ww'$ implies $\ww \eqpm \ww'$ and
 $\ww\ii \rvr {\ww'}\ii $.
\end{lemm}

\begin{proof}
 It suffices to prove the result for $\ww \rvr^{(1)} \ww'$. 
 The case when some factor~$\ss\ii \ss$ has been
 deleted is obvious. Otherwise, assume that $\ww'$ 
 has been obtained from~$\ww$ by
 substituting $\ss\ii \tt$ with $v u\ii $ where $\ss v
 = \tt v$ is one of the relations of the considered
 presentation. Then we have $\ss v \eqp \tt u$, and,
 a fortiori, $\ss v \eqpm \tt u$, hence $\ss\ii \tt \eqpm v
 u\ii $, and, therefore, $\ww \eqpm \ww'$. On the other hand, 
 ${\ww'}\ii $ is obtained from~$\ww\ii $ by replacing~$\tt\ii 
 \ss$ with~$u v\ii $, which is also an $r$-reversing.
\end{proof}

\begin{figure}[htb]
 \includegraphics{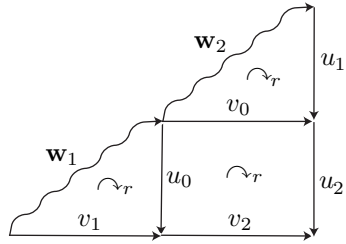}
 \caption{Redressing a product (general case)}
 \label{F:conf}
\end{figure}

\begin{lemm}\label{L:deco}
 (Fig.~\ref{F:conf}) Assume $\ww \rvr^{(k)} v u\ii $
 with $u$, $v \in \SS^*$.
 Then, for every decomposition $\ww = \ww_1
\ww_2$, there exist in~$\SS^*$ decompositions
 $u = u_1 u_2$, $v = v_1 v_2$, and $u_0$, $v_0$ 
 satisfying
 $\ww_1 \rvr^{(k_1)} v_1 u_0\ii $, $\ww_2 \rvr^{(k_2)} v_0
 u_1\ii $, and $u_0\ii v_0 \rvr^{(k_0)} v_2 u_2\ii $ with $k = 
 k_1 + k_2 + k_0$.
\end{lemm}
 
\begin{proof}
 We use induction on~$k$. For $k = 0$, the only
 possibility is $\ww \in \SS^*$, in which case
 we have $\ww = v$ and $u = \e$, and the result is trivial, 
 or $\ww \in (\SS^*)\ii $, in which case we have $\ww
 = u\ii $ and $v = \e$, and the result is trivial as well.
 For $k = 1$, we must have $\ww = \ss\ii \tt$ for some
 letters~$\ss$, $\tt$ such that $\ss v = \tt v$ belongs
 to~$\RR$ (or we have $\ss = \tt$), and everything is
 clear: the result is trivial if either $\ww_1$ or $\ww_2$ is
 empty, and, for $\ww_1 = \ss\ii $ and $\ww_2 = \tt$, we can
 take $u_0 = \ss$ , $v_0 = \tt$, $u_1 = v_1 = \e$, $u_2 = u$,
 and $v_2 = v$, corresponding to $k_0 = 1$, $k_1 = k_2 =
 0$. Assume now $k \ge 2$, and let $\ww'$ be the
 second word in a shortest reversing sequence from~$\ww$
 to~$v u\ii $: by definition, we have $\ww = \uu \ss\ii \tt
 \vv$
 and $\ww' = \uu v' {u'}\ii \vv$, with $\ss$, $\tt $ in~$\SS$
 and $\ss v' = \tt u'$ in~$\RR$. Let us consider a
 decomposition~$\ww =
 \ww_1 \ww_2$. Three cases may happen.
 
 If $\uu \ss\ii \tt$ is a prefix of~$\ww_1$, say $\ww_1 =
 \uu
 \ss\ii \tt \vv_1$, then we have $\ww_1 \rvr^{(1)}\ww'_1$
 with $\ww'_1 = \uu v' {u'}\ii \vv_1$. By construction, we
 have $\ww' = \ww'_1 \ww_2$. Applying the induction
 hypothesis to~$\ww' \rvr^{(k-1)} v u\ii $, we find $u_0$, \pp,
 $v_2$ satisfying $u = u_1 u_2$, $v = v_1 v_2$, and 
 $\ww'_1 \rvr^{(k'_1)} v_1 u_0\ii $, $\ww_2 \rvr^{(k_2)} v_0
 u_1\ii $, $u_0\ii v_0 \rvr^{(k_0)} v_2 u_2\ii $ with $k'_1
 + k_2 + k_0 = k-1$. Then $\ww_1 \rvr^{(1)} \ww'_1$ implies
 $\ww_1 \rvr^{(k'_1 + 1)} v_1 u_0\ii $, and we 
 are done.
 
\begin{figure}
 \includegraphics{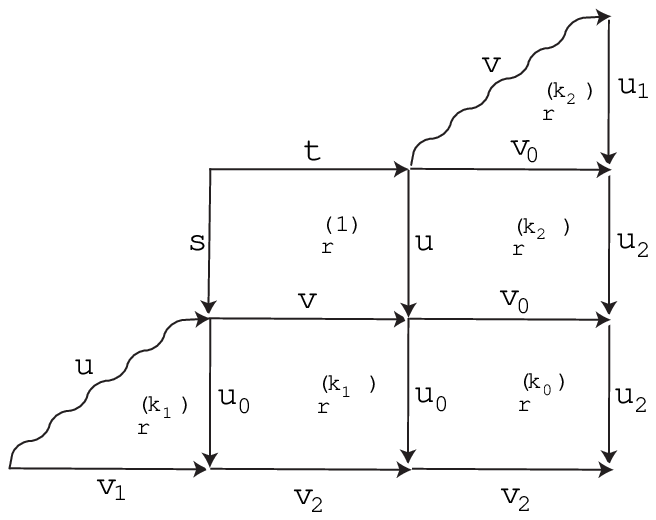}
 \caption{Proof of Lemma~\ref{L:deco}}\label{F:deco}
\end{figure}

 The case when $\ss\ii \tt \vv$ is a suffix
of~$\ww_2$ is 
 symmetric. So we are left with the nontrivial case, namely
 $\ww_1 = \uu \ss\ii $ and $\ww_2 = \tt \vv$
 (Fig.~\ref{F:deco}). Applying the induction hypothesis to
 $\ww' \rvr^{(k-1)} v u\ii $ gives us words $u'_0$,
 $u'_1$, $u'_2$, $v'_0$, $v'_1$, $v'_2$ in~$\SS^*$ satisfying
 $u = u'_1 u'_2$, $v = v'_1 v'_2$, $\ww'_1 \rvr^{(k'_1)}
 v'_1 {u'_0}\ii $, $\ww'_2 \rvr^{(k'_2)} v'_0 {u'_1}\ii $,
 ${u'_0}\ii v'_0 \rvr^{(k'_0)} v_2 u_2\ii $ with $k'_1 + k'_2
 + k'_0 = k-1$. Now, applying the induction hypothesis to 
 $\ww'_1 \rvr^{(k'_1)} v'_1 {u'_0}\ii $ with the decomposition
 $\ww'_1 = \uu v'$ gives us words $u''_0$, $v_1$,
 $v''_2$ in~$\SS^*$ satisfying $v'_1 = v_1 v''_2$, $\uu
 \rvr^{(k''_1)}
 v_1 {u''_0}\ii $, and ${u''_0}\ii v' \rvr^{(k'''_1)} v''_2
 {u'_0}\ii $ with $k''_1 + k'''_1 = k'_1$: indeed, the
 hypothesis that $v'$ belongs to~$\SS^*$ implies that $v'
 \rvr^{(0)} v' \e\ii $ is the only possible reversing from~$v'$.
 Similarly, applying the induction hypothesis to 
 $\ww'_2 \rvr^{(k'_1)} v'_0 {u'_1}\ii $ with the decomposition
 $\ww'_2 = {u'}\ii \vv$ gives us words $v''_0$, 
 $u_1$, $u''_2$ in~$\SS^*$ satisfying $u'_1 = u_1 u''_2$,
 $\vv \rvr^{(k''_2)} v_0 u_1\ii $, and ${u'}\ii v''_0
 \rvr^{(k'''_2)} v'_0 {u''_2}\ii $ with $k''_2 + k'''_2 = k'_2$.
 Put $u_0 = \ss u''_0$, $u_2 = u''_2 u'_2$, $v_0 = \tt v''_0$,
 and $v_2 = v''_2 v'_2$. By construction, we have $u = u_1
 u_2$ and $v = v_1 v_2$. Then we find $\ww_1
 = \uu \ss\ii \rvr^{(k''_1)} v_1 u_0\ii $, and
 $\ww_2 = \tt \vv \rvr^{(k''_2)} v_0 u_1\ii $. Finally,
 we obtain
 \begin{align*}
 u_0\ii v 
 = {u''_0}\ii \ss\ii \tt v''_0
 &\rvr^{(1)} {u''_0}\ii v' {u'}\ii v''_0 \\
 &\rvr^{(k'''_1)} v''_2 {u'_0}\ii {u'}\ii v''_0 \\
 &\rvr^{(k'''_2)} v''_2 {u'_0}\ii v'_0 {u''_2}\ii 
 \rvr^{(k'_0)} v''_2 v'_2 {u'_2}\ii {u''_2}\ii 
 = v_2 u_2\ii ,
 \end{align*}
 hence $u_0\ii v \rvr^{(k_0)} v_2 u_2\ii $ with $k_0 = 
 1 + k'''_1 + k'''_2 + k'_0$. As we check $k''_1 + k''_2 + k_0 
 = k$, we are done.
\end{proof}

\begin{figure}[htb]
 \includegraphics{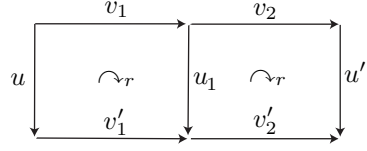}
 \caption{Redressing a product (positive case)}\label{F:prod}
\end{figure}

Applying the previous result to the case when $\ww_1$ has
the form $u\ii v_1$ and $\ww_2$ belongs to~$\SS^*$
gives:

\begin{lemm}\label{L:iter}
 (Fig.~\ref{F:prod}) Assume $u$, $v_1$, $v_2$, $u'$, $v' \in
 \SS^*$ and $u\ii v_1 v_2 \rvr v' {u'}\ii $. Then there
 exists in~$\SS^*$ a decomposition $v' = v'_1 v'_2$ and a
 word~$u_1$ satisfying $u\ii v_1 \rvr v'_1
 {u_1}\ii $ and ${u_1}\ii v_2 \rvr v'_2 {u'}$.
\end{lemm}

\begin{prop}\label{P:eqpo}
 Assume that $(\SS, \RR)$ is a positive group presentation, and $u$,
 $v$, $u'$, $v'$ are words in~$\SS^*$. Then
 $u\ii v \rvr v' {u'}\ii $ implies $u v' \eqp v u'$.
\end{prop}
 
\begin{proof}
 We use induction on the number of steps~$k$ needed
 to reverse $u\ii v$ into~$v' {u'}\ii $. For $k = 0$, the
 only possibility is that $u$ or $v$ is empty, in which case we
 have $u' = u$ and $v' = v$, and the result is true.
 For $k = 1$, the only possibility is that $u$ and $v$
 have length~$1$, \ie, they are letters, say $\ss$, $\tt$
 respectively. In this case, for $\ss\ii \tt$ to reverse to
 $v' {u'}\ii $ means that $\ss v' = \tt u'$ is a relation of the
 presentation, and 
 $\ss v' \eqp \tt u'$ holds by definition.
 Assume now $k \ge 2$. At least one of~$u$, $v$
has
 length larger than~$1$. Assume for instance $\lg(v) \ge 2$,
 and consider a decomposition $v = v_1 v_2$ with $\lg(v_i)
 < \lg(v)$. Applying Lemma~\ref{L:iter}, we obtain $u'_1$,
 $v'_1$, $v'_2$ satisfying $v' = v'_1 v'_2$, $u\ii v_1 \rvr v'_1
 {u_1}\ii $, and ${u_1}\ii v_2 \rvr v'_2 {u'}$
 (Fig.~\ref{F:prod}). The induction
 hypothesis applies to the previous relations, and it gives
 $u v'_1 \eqp v_1 u_1$ and $u_1 v'_2 \eqp v_2 u'$, hence
 $u v'_1 v'_2 \eqp v_1 u_1 v'_2 \eqp v_1 v_2 u'$, \ie, 
 $u v' \eqp v u'$.
\end{proof}

For future use, let us state two applications of the
previous result:

\begin{lemm}\label{L:weak}
 (i) The relation $u\ii w w\ii 
 v \rvr v' {u'}\ii $ implies $u v' \eqp v u'$.
 
 \noindent (ii) The relation $(u v')\ii (v u') \rvr \e$
 implies that there
 exist $u''$, $v''$, $w''$ in~$\SS^*$ satisfying $u\ii v
 \rvr v'' {u''}\ii $, $u' \eqp u'' w''$, and $v' \eqp v''
 w''$.
\end{lemm}

\begin{figure}[htb]
 \includegraphics{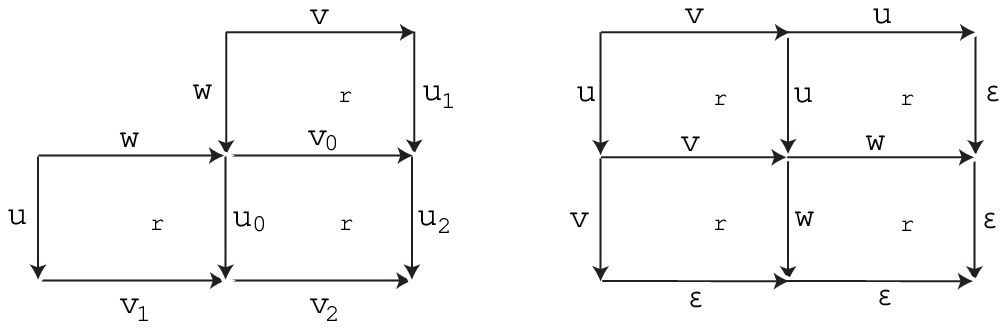}
 \caption{Proof of Lemma~\ref{L:weak}}
 \label{F:lemm}
\end{figure}

\begin{proof}
 (Fig.~\ref{F:lemm})
 (i) Using Lemma~\ref{L:iter}, we see that
 $u\ii w w\ii v \rvr v' {u'}\ii $ implies the existence
 of two decompositions $u' = u'_1 u'_2$, $v' = v'_1
 v'_2$ and of $u_0$, $v_0$ satisfying 
 $u\ii w \rvr v'_1 u_0\ii$, $w\ii v \rvr v_0 {u'_1}\ii$, 
 and $u_0\ii v_0 \rvr v'_2 {u'_2}\ii$.
 Then, using Prop.\ref{P:eqpo}, we obtain
 \[
 u v' = u v'_1 v'_2 \eqp w u_0 v'_2
 \eqp w v_0 u'_2 \eqp v u'_1 u'_2 = v u'.
 \]
 
 (ii) Using Lemma~\ref{L:iter} again, we see that
 $u\ii w w\ii v \rvr v' {u'}\ii $ implies the existence
 of words $u''$, $v''$, $w'$, $w''$ satisfying $u\ii v
 \rvr v'' {u''}\ii $, ${v'}\ii v'' \rvr w'$, ${u''}\ii u' \rvr
 w''$, and ${w'}\ii w'' \rvr \e$. By
 Prop.~\ref{P:eqpo}, the latter relations imply
 $u' \eqp u'' w''$, $w' \eqp w''$, and 
 $v' \eqp v'' w' \eqp v'' w''$.
\end{proof}

The question of whether reversing converges, \ie, the
existence of an upper bound for the length of the
reversing sequences starting from a given word, is
difficult in general. It is easy to give examples of
simple finite presentations, such as the
Baumslag-Solitar presentation $(a, b \; ba = a^2b)$, or
the non-spherical Artin presentation $(a, b, c \;
aba=bab, bcb=cbc, aca=cac)$, where infinitely long
reversing sequences exist: start for instance with $b\ii
ab$ and with $a\ii bc$ in the examples above. Also, 
\cite{Dgd} contains an example of an infinite
presentation where all reversing sequences are finite,
but the only known bound on the length of a
reversing sequence starting from a length~$n$ word
is a tower of exponentials of height~$O(2^n)$.
Besides such complicated cases, easy upper bounds
can be established when the closure of the initial
alphabet under reversing happens to be known.

\begin{defi}
 Assume that $(\SS, \RR)$ is a positive group presentation.
 We say that a subset~$\SSp$ of~$\SS^*$ is {\it closed
 under~$r$-reversing} if $u'$ and
 $v'$ lie in~$\SSp$ whenever $u$ and $v$ do and
 $u\ii v \rvr v' {u'}\ii $ holds. The {\it closure}
 of~$\SS$ under $r$-reversing is defined to be
 the smallest subset of~$\SS^*$ that includes~$\SS$
 and is closed under $r$-reversing.
\end{defi}

\begin{exam}\label{X:clos}
 Let us consider the presentation $(a, b \; a^2 = b^2, 
 ab = ba)$ of Example~\ref{X:aabb}. Then the set
 $\{\e, a, b\}$ is closed under $r$-reversing: up
 to a symmetry, the only possibilities are $a\ii a \rvr
 \e$, $a\ii b \rvr b a\ii$, $a\ii b \rvr a b\ii$, and the
 only words of~$\{a, b\}^*$ involved in the right
 hand sides are $\e$, $a$, and~$b$.
 So $\{\e, a, b\}$ is the closure of~$\{a, b\}$ under
 $r$-reversing.
\end{exam}

Starting with a finite (or, simply, recursive) positive 
group presentation~$(\SS, \RR)$, determining the
closure of~$\SS$ under $r$-reversing is typically a
recursively enumerable process: for each word~$\ww$
on~$\SS \cup \SS\ii$, we can enumerate all
words~$\ww'$ to which $\ww$ is $r$-reversible in
$1$, $2$, etc. steps, and, each time we find a word of
the form~$u v\ii$ with $u$, $v$ in~$\SS^*$, add
it to the current family. Provided we enumerate
the words in a systematic way, all words in the closure
of~$\SS$ will appear at some finite step of the
process, but, if we have no recursive
upper bound for the lengths of the $r$-reversing
sequences from~$\ww$ in terms of the length
of~$\ww$, we shall never know whether all words in
the closure of~$\SS$ have been
found (even if the latter is finite). However, if we
happen to find a finite set of words~$\SSp$ that
includes~$\SS$ and we can prove that every
$r$-reversing sequence from~$u\ii v$ with
$u$, $v \in \SSp$ either ends up with a failure, \ie, with
a word containing some factor~$\ss\ii \tt$ for which
there is no relation $\ss \cdots = \tt \cdots$ in~$\RR$,
or with a word~$v' {u'}\ii$ with $u'$, $v' \in \SSp$,
then we can claim that $\SSp$ includes the closure
of~$\SS$ under reversing. Example~\ref{X:clos}
provides a (trivial) instance of this situation.

\begin{prop}\label{P:dewp}
 Assume that $(\SS, \RR)$ is a recursive positive 
 presentation such that the closure~$\SSr$ of~$\SS$
under $r$-reversing and the
 restriction~$\rvr^{\!\bullet}$ of the relation $u\ii v \rvr v'
 {u'}\ii$ to~${\SSr}^4$ are recursive. Then the
 relation $\ww \rvr v u\ii$ on $(\SS \cup \SS\ii)
 \times (\SS^*)^2$ is recursive; if $\ww$
 is a word with $p$~letters in~$\SS$ and
 $q$~letters in~$\SS\ii$, and $\ww \rvr v u\ii$ holds,
 then $u$ belongs to~${\SSr}^p$, $v$ belongs
 to~${\SSr}^q$, and the reversing of~$\ww$ to~$v
 u\ii$ can be decomposed into at most~$pq$
 reversings in~$\rvr^{\!\bullet}$.
\end{prop}

\begin{proof}
 By hypothesis, the word~$\ww$ is $w_1^{e_1}
 \cdots w_n^{e_n}$ with $w_i \in \SSr$ and $e_i =
 \pm 1$ for $i = 1, \pp, n$. Denote by $d(\ww)$ the
 number of pairs~$(i, j)$ with $i <j$, $e_i = -1$,
 and $e_j = 1$. By construction, we have $d(\ww) \le
 pq$. We prove the result using induction
 on~$d(\ww)$. For $d(\ww) = 0$, the word~$\ww$
 has the form~$v u\ii$ with $v \in {\SSr}^q$
 and $u \in {\SSr}^p$, and it is reversed, so the
 result is true. Otherwise, there exist~$i$ satisfying
 $e_i = -1$ and $e_{i+1} = +1$. Using
 Lemma~\ref{L:deco} twice, we see that there must
 exist $w'_i$, $w'_{i+1}$ in~$\SS^*$ such that
 $w_i\ii w_{i+1} \rvr w'_i {w'}_{\!i+1}^{-1}$ holds and
 $\ww \rvr v u\ii$ may be decomposed into
 \[
 \ww = \uu w_i\ii w_{i-1} \vv
 \rvr \uu w'_i {w'}_{\!i+1}^{-1} \vv \rvr v u\ii
 \]
 with $\uu = w_1^{e_1} \cdots w_{i-1}^{e_{i-1}}$
 and $\vv = w_{i+1}^{e_{i+1}} \cdots 
 \cdots w_{n}^{e_{n}}$. By construction, the words
 $w'_i$ and $w'_{i+1}$ belong to~$\SSr$, and,
 letting $\ww' = \uu w'_i {w'}_{\!i+1}^{-1} \vv$, 
 we see that the word~$\ww'$ satisfies the same
 requirements as~$\ww$ with $d(\ww') = d(\ww) -1$,
 so we can apply the induction hypothesis.
\end{proof}

A favourable case is when all relations in the
considered presentation involve words of length~$2$
at most: in this case, the closure~$\SSr$ of~$\SS$
under reversing is merely $\SS \cup \{\e\}$, so we
obtain:

\begin{coro}
 Assume that $(\SS, \RR)$ is a finite positive
 presentation and all relations in~$\RR$ have the form
 $u = v$ with $u$ and $v$ of length~$1$ or~$2$.
 Then every $r$-reversing sequence starting with a
 length~$n$ word has length $n^2/4$ at
 most, and all words in such a sequence have
 length~$n$ at most.
\end{coro}

The case above is not the only one
when the closure can be determined. For instance, in
the case of the standard presentation of the braid
group~$B_n$, the closure of the generators
$\s_1$, \pp, $\s_{n-1}$ under $r$-reversing is the set
of the $(n-1)!-1$ proper divisors of~$\D_n$. We refer
to~\cite{Pic, Pid, PicNote} for many other examples
(in the complemented case).

\begin{rema}
 It is proved in~\cite{Dgc} that, if $(\SS, \RR)$ is a
 finite complemented presentation and all
 relations in~$\RR$ preserve the length, then there
 exists a constant~$C$ such that, if $u$ and $v$ are
 $\eqp$-equivalent length~$n$ words, then $u\ii v$
 reverses to~$\e$ in at most $2^{2^{Cn}}$~steps.
 Whether this result extends to arbitrary
 finite presentations is unknown.
\end{rema}

\section{Complete presentations}\label{S:comp}

We introduce now our key notion, namely that of a
complete presentation. The idea is that a presentation
is complete if it contains enough relations to
make reversing exhaustive.

\begin{defi}\label{D:comp}
 Let $(\SS, \RR)$ be a positive presentation. For $u$, 
 $v$, $u'$, $v'$ in~$\SS^*$, we say that $(\SS, \RR)$
 is {\it $r$-complete} at $u$, $v$, $u'$, $v'$ if the
 following implication holds:
 \begin{equation}\label{E:rcom}
 \begin{matrix}
 \text{If $uv' \eqp vu'$ holds, then there exist
 $u''$, $v''$, $w$ in~$\SS^*$} \hfill\\
 \text{satisfying
 $u\ii v \rvr v'' {u''}\ii$, 
 $u' \eqp u'' w$, 
 and $v' \eqp v'' w$;}
 \end{matrix}
 \end{equation}
 we say that $(\SS, \RR)$ is $r$-complete if
 \eqref{E:rcom} holds for all $u$, $v$, $u'$, $v'$. 
 
 Symmetrically, we say that $(\SS, \RR)$
 is {\it $l$-complete} at $u$, $v$, $u'$, $v'$ if we
 have
 \begin{equation}\label{E:lcom}
  \begin{matrix} 
   \text{If $v'u \eqp u'v$ holds, then there exist
   $u''$, $v''$, $w$ in~$\SS^*$} \hfill\\
   \text{satisfying
   $u v\ii \rvl {v''}\ii u''$, 
   $u' \eqp w u''$, 
   and $v' \eqp w v'$;}
  \end{matrix}
 \end{equation}
 we say that $(\SS, \RR)$ is $l$-complete if 
 \eqref{E:lcom} holds for all $u$, $v$, $u'$, $v'$, and
 that $(\SS, \RR)$ is {\it complete} if it is both $r$- and
 $l$-complete.
\end{defi}

\begin{figure}[htb]
 \includegraphics{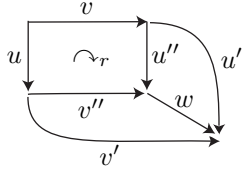}
 \caption{The $r$-completeness condition}
 \label{F:comp}
\end{figure}

Completeness says something nontrivial only for those
$4$-tuples that satisfy $u v' \eqp v u'$: for the other
ones, the implications \eqref{E:rcom} and
\eqref{E:lcom} are trivially true. By
Prop.~\ref{P:eqpo}, $u\ii v \rvr v'' {u''}\ii
$ and, symmetrically, $u v\ii \rvl {v''}\ii u''$ imply $u
v'' \eqp v u''$, so the converse implications
of~\eqref{E:rcom} and \eqref{E:lcom} always hold.
Completeness claims that these sufficient conditions
also are necessary: it tells us that every
common multiple relation~$u v' \eqp v u'$ factors
through some reversing, as illustrated in
Fig.~\ref{F:comp}.

\begin{rema}
 The statement of the completeness property and the
 picture in Fig.~\ref{F:comp} are formally reminiscent
 of Prop.~H in~\cite{Gar}, or of the
 K\"urzungslemma in~\cite{BrS}, or of the chainability
 condition of~\cite{Cor}. However, the point here is
 not the factorization property for common multiples,
 but the fact that the square $(u v', v u')$ corresponds
 to an $r$-reversing process: completeness is a
 property of a presentation, not of a monoid.
\end{rema}

The following result is a straightforward consequence
of the definition:

\begin{lemm}\label{L:complextends}
 Assume that $(\SS, \RR)$ is an $r$-complete positive
 presentation, and $\RR'$ includes~$\RR$. Then $(\SS,
\RR')$ is $r$-complete as well.
\end{lemm}

A natural question is whether complete presentations
exist. The answer is trivial:

\begin{prop}\label{P:exis}
 Every monoid with no nontrivial unit admits a
 complete presentation.
\end{prop}

\begin{proof}
 Let $M$ be a monoid, and $\SS$ be an arbitrary
 set of generators for~$M$. Let~$\cong$ be the
 congruence on~$\SS^*$ such that $M$ is the
 quotient $\SS^* \!/\! \cong$. Let $\RR$ consist of
 all relations $u = v$ with $u \cong v$ and $u$,
 $v \not= \e$. As $u \cong \e$ is supposed to hold
 for no nonempty word~$u$, $(\SS, \RR)$ is a
 presentation of~$M$, which we claim is complete.
 Indeed, assume $u v' \eqp v u'$. If $u$ or $v$ is
 empty, the condition for completeness holds trivially.
 Otherwise, we write $u = \ss u_0$, $v = \tt v_0$ with $\ss$,
 $\tt \in \SS$. The hypothesis is $\ss u_0 v' \eqp
 \tt v_0 u'$, hence the relation $\ss u_0 v' = \tt
 v_0 u'$ belongs to~$\RR$ as the
 considered words are nonempty. Then $\ss\ii \tt \rvr u_0
 v' {u'}\ii v_0\ii$ holds by definition, which implies
 \[
 u\ii v = u_0\ii \ss\ii \tt v_0 
 \rvr u_0\ii u_0 v' {u'}\ii v_0\ii v_0 = v' {u'}\ii .
 \]
 Putting $u'' = u'$, $v'' = v'$ and $w = \e$ gives
 \eqref{E:rcom}, proving $r$-completeness at $u$,
 $v$, $u'$, $v'$.
 The verification of $l$-completeness is similar.
\end{proof}

The practical interest of the previous result is weak: the
complete presentation given by Prop.~\ref{P:exis}
is infinite whenever the considered monoid is infinite,
and, more important, writing such a presentation
supposes knowing a solution to the word problem. As
we shall see below, the interesting case is that of a finite
complete presentation, about which Prop.~\ref{P:exis}
tells us nothing in general.

A more interesting method to possibly obtain complete
presentations consists in considering minimal
common multiples (when they exist).

\begin{defi}
 Assume that $M$ is a monoid. For $\xx$, $\yy$, $\zz
 \in M$, we say that $\zz$ is a {\it minimal} common
 right multiple, or {\it $r$-mcm}, of~$\xx$ and~$\yy$
 if $\zz$ is a right multiple both of~$\xx$ and~$\yy$,
 but no proper left divisor of~$\zz$ is. 
\end{defi}

The notion of a minimal common multiple is a
generalization of that of a least common multiple:
saying that two elements $\xx$, $\yy$ admit a least
common multiple amounts to saying that they admit a
unique minimal common multiple. Mcm's need not
exist in general, but they do in good cases, namely
when the considered monoid is {\it Noetherian}. If
$\xx$, $\yy$ are elements of a monoid~$M$, we write
$\xx \divl \yy$ if $\yy = \xx \zz$ holds for some~$\zz
\not= 1$, and, symmetrically,
$\xx \divr \yy$ if $\yy = \zz \xx $ holds for some~$\zz
\not= 1$. 

\begin{defi}
 We say that a monoid~$M$ is {\it
 $l$-Noetherian} if the relation~$\divl$ has no
 infinite descending chain, \ie, there exists no infinite
 sequence $\xx_0 \multl \xx_1\multl \pp$ in~$M$.
 Symmetrically, we say that $M$ is {\it $r$-Noetherian}
 if $\divr$ has no infinite descending chain, and that
 $M$ is {\it Noetherian} if it is both $l$- and
 $r$-Noetherian.
\end{defi}

If $M$ is an $l$-Noetherian monoid, the associated
relation~$\divl$ must be irreflexive, so, in particular,
$M$ contains no nontrivial invertible element;
more generally, the relation~$\divl$ is then a partial
ordering on~$M$, which is compatible with
multiplication on the left, and for which $1$ is a least
element. 

\begin{lemm}\label{L:mcmexist}
 Assume that $M$ is an $l$-Noetherian monoid. Then
 any common $r$-multiple of two elements~$\xx$,
 $\yy$ of~$M$ is an $r$-multiple of some $r$-mcm
 of~$\xx$ and~$\yy$.
\end{lemm}

\begin{proof}
 Our hypothesis is that every nonempty subset of~$M$
 contains an element which is minimal with respect
 to~$\divl$. Applying this property to the set of all
 common right multiples of~$\xx$ and~$\yy$ which
 are left divisors of~$\zz$ gives the expected right
 mcm.
\end{proof}

We shall now prove how to obtain complete
presentations in the case of a Noetherian monoid by
considering $r$-mcm relations.

\begin{defi}
 Assume that $M$ is a monoid, and $\SS$ is a set of
 generators for~$M$. We say that a family of
 relations~$\RR$ is an {\it $r$-selector} on~$\SS$
 in~$M$ if, for all~$\ss$,~$\tt$ in~$\SS$ and for
 each $r$-mcm~$\xx$ of~$\ss$ and~$\tt$, there
 exists one pair of words~$(u, v)$ in~$\SS^*$
 such that $\ss v = \tt u$ belongs to~$\RR$ and both
 $\ss v$ and $\tt u$ represent~$\xx$.
\end{defi}

Thus, an $r$-selector is a family of relations that proves
all equalities connected with right mcm's in the
considered monoid~$M$. Observe that $r$-selectors
always exist, but an $r$-selector may be just empty
when no right mcm exists. The following result shows
that, in the case of a Noetherian monoid, each
$r$-selector gives rise to a presentation, which
moreover turns out to be $r$-complete. 

\begin{prop}\label{P:exno}
 Assume that $M$ is a left cancellative Noetherian
 monoid, $\SS$ is a set of generators for~$M$, and
 $\RR$ is an $r$-selector on~$\SS$ in~$M$. Then $(\SS, \RR)$ is an $r$-complete presentation of~$M$.
\end{prop}

\begin{proof}
 As in the proof of Prop.~\ref{P:exis}, let $\cong$
 denote the congruence on~$\SS^*$ such that $M$ is
 isomorphic to~$\SS^* \! / \! \cong$. Let $\eqp$
 be the congruence associated with the selector~$\RR$. By
 definition, $\RR$ consists of pairs $\{u, v\}$ that satisfy
 $u \cong v$, so $u \eqp v$ implies $u \cong v$
 trivially. 
 
 We shall now prove conversely that $u \cong v$
 implies $u \eqp  v$ for all $u$, $v$ in~$\SS^*$ using
 induction
 on~$\cl u$ with respect to~$\divr$ (we recall that
 $\cl u$ denotes the element of~$M$ represented
 by~$u$). As $1$ is the least element relative
 to~$\divr$, let us first assume $\cl u = \cl v = 1$, \ie,
 $u \cong v \cong \e$. We have seen that $1$ is the
 only invertible element in~$M$, so, necessarily, $u$
 and $v$ are empty, and we have $u = v$, hence $u
 \eqp v$.

 Assume now $\cl u = \cl v \multr 1$. Then $u$
 and $v$ are nonempty words, say $u = \tt u_0$, $v =
\ss v_0$, with $\ss$, $\tt \in \SS$. The hypothesis $u
\cong v$ means that $\cl u$ is a common
$r$-multiple of~$\ss$ and~$\tt$ in~$M$. By
Lemma~\ref{L:mcmexist}, some left divisor~$\zz$
of~$\cl u$ has to be an $r$-mcm of~$\ss$ and~$\tt$.
So, by definition, there must exist some relation $\tt u'
= \ss v'$ in~$\RR$ such that both $\tt u'$ and $\ss v'$
represent~$\zz$ in~$M$, and the hypothesis that
$\zz$ is a left divisor of~$\cl u$ implies that  some
word~$w$ satisfies
$$\tt u_0 \cong \tt u' w \cong \ss v' w \cong \ss
v_0.$$
Applying the hypothesis that $M$ is left cancellative,
we deduce $u_0 \cong u' w$ and $v_0 \cong v' w$.
By construction, $\cl{u_0}$ and $\cl{v_0}$ are proper
right divisors of~$\cl u$, so the induction hypothesis
allows us to deduce $u_0 \eqp u' w$ from $u_0
\cong u' w$ and $v_0 \eqp v' w$ from $v_0 \cong
v' w$, and we obtain
$$u = \tt u_0 \eqp \tt u' w \eqp \ss v' w \eqp \ss
v_0 = v.$$
 
 It remains to prove that $(\SS, \RR)$ is
 $r$-complete: we postpone the proof to
 Sec.~\ref{S:crit}, as the needed argument is similar to,
 but simpler than, the argument developed for
 Prop.~\ref{P:loct}.
\end{proof}

In order to obtain a complete presentation (and not
only an $r$-complete one), we can appeal to the
symmetric obvious notion of an $l$-selector, and
using Proposition~\ref{P:exno}, its left counterpart, and
Lemma~\ref{L:complextends}, we obtain

\begin{prop}
 Assume that $M$ is a cancellative Noetherian monoid,
 $\SS$ is a set of generators for~$M$, $\RR_r$ is an
 $r$-selector on~$\SS$ in~$M$, and $\RR_l$ is an
 $l$-selector on~$\SS$ in~$M$. Then $(\SS, \RR_r \cup \RR_l)$ is a complete presentation of~$M$.
\end{prop}

Let us conclude this section with yet another way of
constructing a complete presentation, even in a
non-Noetherian case, when what is called a spanning
subset in~\cite{Dgo} happens to be known.

\begin{prop}\label{P:span}
 Assume that $M$ is a monoid, $\SS$ is
 a set of generators for~$M$, and $\SS'$ is a
 subset of~$M$ that includes~$\SS$ and satisfies the
 following condition: 
 
 {\narrower\noindent $(*)$ For all $\xx$, $\yy$
 in~$\SS'$, if $\zz$ is a common right multiple
 of~$\xx$ and~$\yy$, then there exist $\xx'$ and
 $\yy'$ in~$\SS'$ satisfying $\xx \yy' = \yy \xx' \divel
 \zz$. \par}
 
 \noindent Let $\RR'$ be the set of
 all relations $\xx \yy' = \yy \xx'$ and $\xx \yy' = \yy$
 with $\xx$, $\yy$, $\xx'$, $\yy' \in \SS'$. Then, 
 for each~$\xx$ in~$\SS'$, let $f(\xx)$ be a word
 in~$\SS^*$ representing~$\xx$, and let~$\RR$ be
 the  image of~$\RR'$ under~$f$. Then $(\SS', \RR')$
 and $(\SS, \RR)$ are $r$-complete presentations
 of~$M$.
\end{prop}

\begin{proof}
 That $(\SS', \RR')$ is a presentation of~$M$ is proved
 in~\cite{Dgo}. The argument is similar to that of
 Prop.~\ref{P:exno}, but it uses an induction on
 $\lg(u) + \lg(v)$ instead of an induction
 on~$\cl u$ relative to~$\divr$, which need not be
 well-founded. The
 $r$-completeness of the presentation is then a
 direct translation of the hypothesis on~$\SS'$.
 
 As for~$(\SS, \RR)$, by construction, every relation
 in~$\RR'$ follows from one relation in~$\RR$, so $(\SS,
 \RR)$ is a presentation of~$M$. Assume
 $u v' \eqp^{\RR} v u'$. As $\SS$ is included 
 in~$\SS'$, the words~$u$, $v$, $u'$, $v'$ are words
 on~$\SS'$, and we have $u v' \eqp^{\RR'} v u'$
 since $(\SS', \RR')$ is a presentation of~$M$.
 As $(\SS', \RR')$
 is $r$-complete, we must have $u\ii v \rvr^{\RR'} v'
 {u'}\ii$, $u' \eqp^{\RR'} u'' w$, and $v' \eqp^{\RR'} v'
  w$ for some words $u''$, $v''$, $w$ in~${\SS'}^*$.
 Then $u\ii v \rvr^{\RR'} v' {u'}\ii$ implies $u\ii v 
 \rvr^{\RR} f(v') f(u')\ii$, as we observe as in the proof
 of Lemma~\ref{L:rest} that, if $\ss$ and $\tt$ are two
  letters in~$\SS'$ and ${\ss}\ii \tt \rvr^{\RR'} v u\ii$
 holds, then ${f(\ss)}\ii f(\tt) \rvr^{\RR} f(v) f(u)\ii$
 holds as well, and then use an induction on the
 number of reversing steps. Next $u' \eqp^{\RR'} u''
 w$ implies $u\ii v \rvr^{\RR} f(v')
 f(u')\ii$ by definition of~$\RR$, and, similarly, 
 $v' \eqp^{\RR'} v' w$ implies $v' \eqp^{\RR}
 f(v'') f(w)$. This shows that the words $f(u'')$, $f(v'')$,
 and $f(w)$ fulfill the requirements for $(\SS, \RR)$ to
 be $r$-complete at $u$, $v$, $u'$, $v'$.
\end{proof}

As the connection between Condition~$(*)$ in
Prop.~\ref{P:span} and $r$-completeness is
clear, the previous result is essentially trivial, and
so is the converse statement that, if $(\SS, \RR)$ is an
$r$-complete presentation of some mon\-oid~$M$,
then the subset of~$M$ consisting of those elements
that can be represented by words in the closure
of~$\SS$ under $r$-reversing satisfies
Condition~$(*)$.

\section{The cube condition}\label{S:cube}

At this point, we know that every monoid~$M$ 
(with no nontrivial unit) and every group~$G$ admit
complete presentations, but we are left with the
question of recognizing that a given presentation is
possibly complete. In every case, even for a finite
presentation, the question is nontrivial, as checking
$r$-completeness for one particular 4-tuple of words
requires being able to decide $\eqp$-equivalence, and
checking it for all 4-tuples is an infinite process.

In this section, we introduce a new combinatorial
condition involving reversing, the cube condition,
and we prove that completeness is equivalent to
that cube condition being satisfied. This is
a first step toward an effective completeness criterion
that will be established in the subsequent section.

\begin{defi}\label{D:cube}
 Let $(\SS, \RR)$ be a positive presentation. For $u$, 
 $v$, $w$ in~$\SS^*$, we say that $(\SS, \RR)$ satisfies
 the {\it $r$-cube condition} (\resp the {\it strong
 $r$-cube condition}) at~$u$, $v$, $w$ if the
 implication
 \begin{equation}\label{E:rcub}
  \begin{matrix}
   \text{If we have $u\ii w w\ii v \rvr v' {u'}\ii $ 
   with $u'$, $v' \in \SS^*$,} \hfill\\
   \text{then there exist $u''$, $v''$, $w''$
   in~$\SS^*$ satisfying}\hfill\\
   \qquad \text{$u\ii v \rvr v'' {u''}\ii $, $u' \eqp u''
   w''$, and $v' \eqp v'' w''$.}\hfill\\
   \text{(\resp then we have $(u v')\ii (v u') \rvr \e$);}
   \hfill
  \end{matrix}
 \end{equation}
 for $\SSp \ince \SS^*$, we say that $(\SS,
\RR)$ satisfies the (strong) $r$-cube condition
 {\it on}~$\SSp$ if the (strong) $r$-cube
 condition holds for all $u$, $v$, $w$ in~$\SSp$.
\end{defi}

\begin{figure}[htb]
 \includegraphics{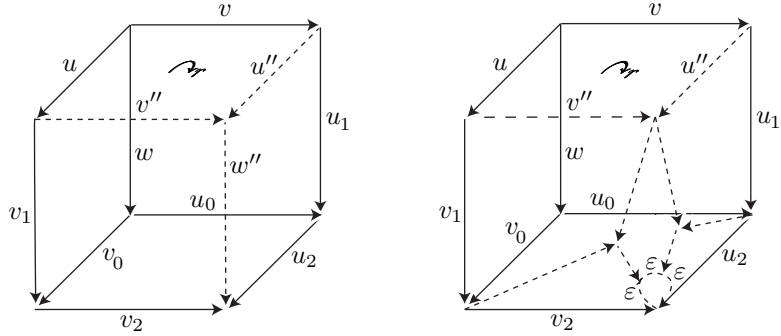}
 \caption{The $r$-cube and the 
 strong $r$-cube conditions}\label{F:cube}
\end{figure}

The cube conditions are illustrated in
Fig.~\ref{F:cube}. We start with an incomplete cube
consisting of three faces constructed on~$(u, w)$,
$(w, v)$, and $(u_0, v_0)$ and correspond to
$r$-reversings, and the condition means that we can
complete the cube with a top reversing face and a
last edge. In the cube condition, we require that the
last two faces correspond to equivalences, while, in
the strong cube condition, we require that the last two
faces correspond, in a slightly more complicated way,
to reversings. As the name suggests, the strong cube
condition implies the cube condition: indeed, 
Lemma~\ref{L:weak} tells us that $(u v')\ii (v u') \rvr 
\e$ implies the existence of $u''$, $v''$, $w''$
satisfying $u\ii v \rvr v'' {u''}\ii$, $u' \eqp u'' w''$, and
$v' \eqp v'' w''$. We shall see below that both
conditions actually are equivalent in the case of an $r$-complete presentation.

\begin{figure}
 \includegraphics{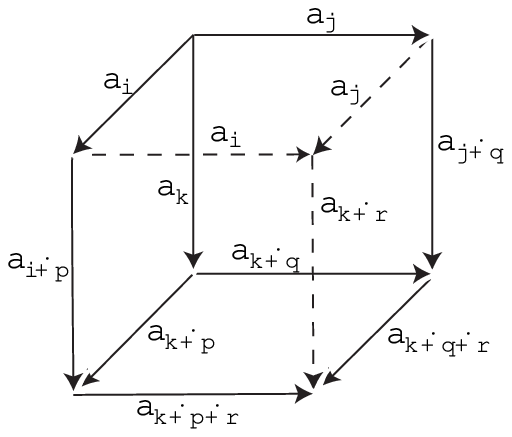}
 \caption{The cube condition for $\MO(\SS_n;
  \RR_n)$}
 \label{F:exam}
\end{figure}

\begin{exam}\label{X:aaii}
 Let $\SS_n$ be $\{a_1, \pp, a_n\}$, and $\RR_n$
 be the family of all relations
 $a_ia_{i \op p} = a_ja_{j \op p}$ with
 $1 \le i < j \le n$ and $1 \le p \le n$, where $x \op
 y$ denotes the unique
 number in $(1, \pp, n)$ equal to
 $x + y$ modulo~$n$. For instance, the
 monoid $\MO(\SS_2; \RR_2)$ is (isomorphic
 to) $\MO(a, b; a^2 = b^2, ab= ba)$
 considered in Example~\ref{X:aabb}, while 
 $\MO(\SS_3; \RR_3)$ is (isomorphic
 to) 
 \[
 \MO(a, b, c; a^2 = b^2 = c^2, 
 ab = bc = ca, 
 ac = ba = cb).
 \]
 We claim that the (strong) $r$-cube condition is
 satisfied
 by~$(\SS_n, \RR_n)$ for every
 triple of letters $a_i$, $a_j$, $a_k$. Indeed, 
 the words to which $a_i\ii a_k$ reverses are the
 words $a_{i \op p} a_{k \op p}\ii $ with $1 \le
 p \le n$; similarly, the words to which
 $a_k\ii a_j$ reverses are the
 words $a_{k \op q} a_{j \op q}\ii $ with $1 \le
 q \le n$; finally, the words to which $a_{k
 \op p}\ii a_{k \op q}$ reverses are the
 words $a_{k \op p \op r} a_{k \op q \op r}\ii $
 with $1 \le r \le n$. But, then, $a_i\ii a_j$ 
 reverses to~$a_i a_j\ii $, and we have
 $a_{i \op p} a_{k \op p \op r} \eqp
 a_i a_{k \op r}$ and 
 $a_{j \op q} a_{k \op q \op r} \eqp
 a_j a_{k \op r}$ (Fig.~\ref{F:exam}), which is the
 $r$-cube condition at $a_i$, $a_j$, $a_k$.
 Moreover, we find
 \begin{gather*}
 a_{k \op p \op r}\ii a_{i \op p}\ii a_i a_{k \op r}
 \rvr a_{k \op p \op r}\ii a_{k \op p \op r}
 a_{k \op r}\ii a_{k \op r} \rvr \e, \\
 a_{k \op q \op r}\ii a_{j \op q}\ii a_j a_{k \op r}
 \rvr a_{k \op q \op r}\ii a_{k \op q \op r}
 a_{k \op r}\ii a_{k \op r} \rvr \e,
 \end{gather*}
 which gives the strong $r$-cube condition.
\end{exam}

The connection between completeness and cube
condition is as follows:

\begin{prop}\label{P:cocu}
 A positive presentation $(\SS, \RR)$ is $r$-complete if
 and only if any of the following four equivalent
 conditions is satisfied:
 
 (i) Equivalence is detected by $r$-reversing: $u \eqp
 v$ is equivalent to $u\ii v \rvr \e$.

 (ii) The relation $u\ii v \rvr \e$ is transitive.

 (iii) The strong $r$-cube condition is satisfied on~$\SS^*$.

 (iv) The $r$-cube condition is satisfied on~$\SS^*$.
\end{prop}

\begin{proof}
 Assume $u \eqp v$, \ie, $u \e \eqp v \e$. If $(\SS,
 \RR)$ is $r$-complete, we obtain
 $u''$, $v''$, $w$ satisfying $u\ii v \rvr v'' {u''}\ii$,
 $\e \eqp u'' w$, and $\e \eqp v'' w$. As $(\SS, \RR)$
 is positive, $\e \eqp u'' w$ implies
 $u'' = w = \e$, and $\e \eqp v'' w$ implies $v'' = \e$.
 This means that we have $u\ii v \rvr \e$, and (i) is
 true.
 
 Conversely, assume~$u v' \eqp v u'$. If (i) holds, we
 have $(u v')\ii (v u') \rvr \e$. By
 Lemma~\ref{L:weak}(ii), we obtain
 $u''$, $v''$, $w'$, and $w''$ satisfying $u\ii v \rvr v''
 {u''}\ii $, $u' \eqp u'' w''$ and $v' \eqp v'' w''$,
 \ie, the $r$-completeness condition for $u$, $v$, $u'$,
 $v'$ is satisfied. So $r$-completeness is equivalent
 to~(i).
 
 Next, by definition, the relation~$\eqp$ is an equivalence
 relation, hence it is transitive, so (i) implies~(ii).
 Conversely, by construction, the relation $u\ii v \rvr \e$ is
 always reflexive, symmetric, and compatible with
 multiplication on both sides so, if (ii) holds, the
 relation is a
 congruence on the monoid~$\SS^*$. By
 Prop.~\ref{P:eqpo}, this congruence is included
 in~$\eqp$. On the other hand, it contains all relations
 of~$\RR$, so it includes~$\eqp$, and,
 finally, it coincides with the latter. So (ii) is
 equivalent to~(i).
 
 Assume now $u\ii w w\ii v \rvr v' {u'}\ii $. By
 Lemma~\ref{L:iter}, there exist $u_i$, $v_i$, $i = 0, 1, 2$,
 satisfying $u\ii w \rvr v_1 u_0\ii $, $w\ii v \rvr v_0 w_1\ii $, 
 and $u_0\ii v_0 \rvr v_2 u_2\ii $, and we have $u'= u_1
 u_2$ and $v' = v_1 v_2$ (as in Fig.~\ref{F:cube}). We read
 \[
 u v' = u v_1 v_2 
 \eqp w u_0 v_2 
 \eqp w v_0 u_2 
 \eqp v u_1 u_2 = v u',
 \]
 hence $u v'
 \eqp v u'$. If the presentation is $r$-complete, this implies
 that there exist $u''$, $v''$, $w''$ satisfying $u\ii v \rvr
 v'' {u''}\ii $, $u' \eqp u'' w''$, and
 $v' \eqp v'' w''$, which gives the strong $r$-cube condition.
 So $r$-completeness, hence (ii) as well, implies~(iii),
 hence~(iv) by Lemma~\ref{L:weak}(i).
 
 Finally, assume $u\ii w \rvr \e$ and $w\ii v
 \rvr \e$. As $\e\ii \e \rvr \e$ trivially holds, we deduce
 $u\ii w w\ii v \rvr \e$. If the $r$-cube condition is
 satisfied, we deduce that there exist $u''$, $v''$, $w''$
 satisfying $u\ii v \rvr v' {u''}\ii $, $\e \eqp
 u'' w''$, and $\e \eqp v'' w''$. The latter relations imply
 $u'' = v'' = w'' = \e$, hence $u\ii v \rvr \e$. This
 shows that (iv) implies~(ii), and, therefore, that (ii),
 (iii), and~(iv) are equivalent.
\end{proof}

By Prop.~\ref{P:cocu}, establishing the
possible completeness of a presentation reduces to
establishing the (strong) cube condition for all
triples of words. Observe that, in practice, checking
the strong cube condition is easier than
checking the cube condition, as the former involves
only reversing, while the latter involves the
equivalence relation~$\eqp$ of which we have no
control as long as the presentation is not known to be
complete.

In the complemented case, \ie, when $r$-reversing
is a deterministic process, the cube condition takes
special forms that have been considered in~\cite{Dff}
and~\cite{Dgk}. Indeed,  in this case, there exists for
each pair of words~$u$, $v$ at most one pair of
words~$(u', v')$ satisfying $u\ii v \rvr v' {u'}\ii $. Let
us define $(u \dr v, v \dr u)$ to be the unique such
pair~$(u', v')$ when it exists---by
Lemma~\ref{L:symm}, the symmetry of reversing
makes the definition unambiguous. 

\begin{lemm}\label{L:cohe}
 Assume that $(\SS, \RR)$ is a complemented
 presentation. Then a sufficient
 condition for the $r$-cube (\resp the
 strong $r$-cube) condition to be satisfied at $u$, $v$,
 $w$ is that the relation
 \begin{gather}
 (u \dr v) \dr (u \dr w) 
 \eqp (v \dr u) \dr (v \dr w) \label{E:cohe}\\ 
 (\quad \resp ((u \dr v) \dr (u \dr w)) 
 \dr ((v \dr u) \dr (v \dr w)) = \e\quad)
 \label{E:cohf}
 \end{gather}
 and the relations obtained by permutation of $u$,
 $v$, $w$ are satisfied.
\end{lemm}

\begin{figure}[htb]
 \includegraphics{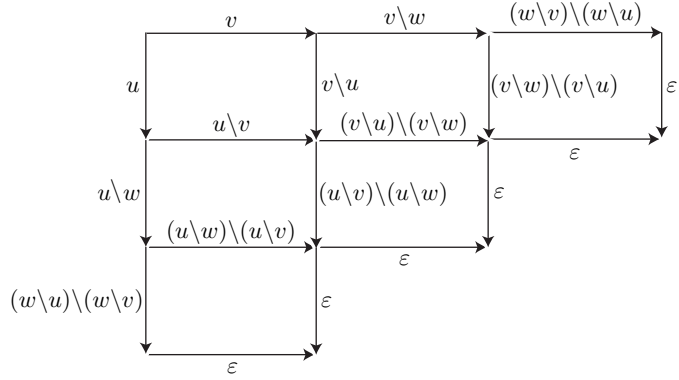}
 \caption{Strong cube condition in the complemented
 case}\label{F:cohe}
\end{figure}

\begin{proof}
 The only word of the form $v' {u'}\ii$
 to which $u\ii w w\ii v$ reverses is 
 \[
 (u\dr w) ((w \dr u) \dr (w \dr v))
 ((w \dr v) \dr (w \dr u))\ii (v \dr w)\ii,
 \]
 and the only word of this form to which $u\ii v$
 may reverse is $(u \dr v) (v \dr u)\ii$. So the point
 for the cube condition is to find~$w'$ satisfying
 \[
 (u\dr w) ((w \dr u) \dr (w \dr v)) \eqp (u \dr v) w'
 \mbox{ and }
 (v\dr w) ((w \dr v) \dr (w \dr u)) \eqp (v \dr u) w'.
 \]
 Now, assuming \eqref{E:cohe} and its cyclic analogs,
 and using the identity $u_1 (u_1 \dr v_1) \eqp
 v_1 (v_1 \dr u_1)$, which is the form taken by
 Prop~\ref{P:eqpo} in this context, we find
 \begin{align*}
 (u\dr w) ((w \dr u) \dr (w \dr v))
 & \eqp (u\dr w) ((u \dr w) \dr (u \dr v)) \\
 & \eqp (u\dr v) ((u \dr v) \dr (u \dr w)) 
 \eqp (u\dr v) ((v \dr u) \dr (v \dr w)), \\
 (v\dr w) ((w \dr v) \dr (w \dr u)) 
 & \eqp (v\dr w) ((v \dr w) \dr (v \dr u)) 
 \eqp (v\dr u) ((v \dr u) \dr (v \dr w)),
 \end{align*}
 the expected form with $w' = (v \dr u) \dr (v \dr
 w)$.
 
 As for the strong cube condition, we wish to prove
 the relation
 \[
 ((w \dr v) \dr (w \dr u))\ii 
 (u \dr w) \ii u\ii v (v \dr w) 
 ((w \dr v) \dr (w \dr u)) \rvr \e.
 \]
 Fig.~\ref{F:cohe} gives the result assuming
 \eqref{E:cohf} and its analogs.
\end{proof}

It is not clear that the sufficient conditions of
Lemma~\ref{L:cohe} are necessary for a given triple
of words~$(u, v, w)$, but they are globally necessary
in that, if \eqref{E:cohe} or \eqref{E:cohf} is satisfied
for all triples~$(u, v, w)$, then, as is proved
in~\cite{Dgk}, $u \eqp v$ is equivalent to $u\ii v \rvr
\e$, so, in our current framework, the presentation is
$r$-complete, and, therefore, the cube and the strong
cube conditions are satisfied for all triples.

\begin{figure}
 \includegraphics{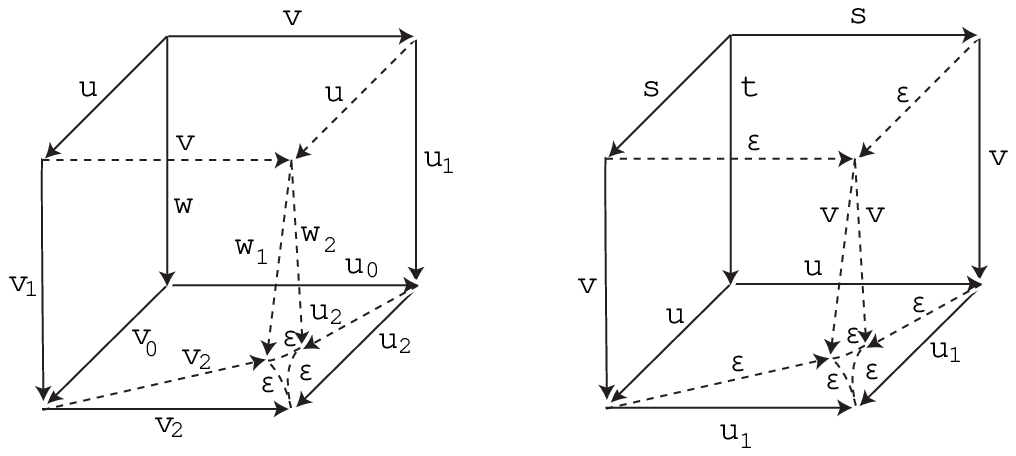}
 \caption{Variants of the cube conditions}
 \label{F:decu}
\end{figure}

\begin{rema}
 The most natural generalization of
 Condition~\eqref{E:cohe} would be: 
 \begin{equation}\label{E:rcuv}
 \begin{matrix}
 \text{Assume $u\ii w \rvr v_1 u_0\ii $, 
 $w\ii v \rvr v_0 u_1\ii $, and $u_0\ii v_0 \rvr v_2 
 u_2\ii $;} \hfill\\
 \text{then there exist $u''$, $v''$, $u''_2$, $v''_2$, 
 $w_1$, $w_2$ in~$\SS^*$ satisfying}\hfill\\
 \qquad \text{$u\ii v \rvr v'' {u''}\ii $, $v_1\ii v'' \rvr 
 v''_2 w_1$, ${u''}\ii u_1 \rvr w_2 {u''_2}\ii,
 $}\hfill\\
 \qquad \text{and $v_2 \eqp v''_2$, $u'_2 \eqp
 u''_2$, $w_1 \eqp w_2$ (Fig.~\ref{F:decu} left).}
 \hfill
 \end{matrix}
 \end{equation}
 However, Condition~\eqref{E:rcuv} is not
 suitable, as it may hold only if the
 considered presentation is equivalent to a
 complemented presentation, at least if there is no
 relation~$\ss \cdots = \ss \cdots$ in~$\RR$ and
 $r$-reversing is convergent, \ie, every word~$u\ii v$
 reverses to at least one word~$v' {u'}\ii $. Indeed, assume
 that $\ss v = \tt u$ and $\ss v' = \tt u'$ belong
 to~$\RR$.
 Then we have $\ss\ii \tt \rvr v u\ii $, $\tt\ii \ss \rvr u'
 {v'}\ii $, and there exist $u_1$, $u'_1$ satisfying $u\ii u' \rvr
 u'_1 u_1\ii $ (Fig.~\ref{F:decu} right). We
apply~\eqref{E:rcuv}: as $\ss\ii \ss \rvr \e$, $v\ii \e
\rvr \e v\ii $, and $\e\ii v'
 \rvr v' \e\ii $ are the only possibilities, and $u_1 \eqp
 \e$ implies $u_1 = \e$, we deduce $v \eqp v'$ and
 $u\ii u' \rvr \e$, hence $u \eqp u'$, \ie, the two
 relations $\ss \cdots = \tt \cdots$ are essentially one
 and the same relation.
 
 The same remark applies to the most natural 
 generalization of Condition~\eqref{E:cohf}, 
 namely the following variant of~\eqref{E:rcuv}
 corresponding to a $6$-face reversing cube:
 \begin{equation}\label{E:rcuw}
 \begin{matrix}
 \text{Assume $u\ii w \rvr v_1 u_0\ii $, 
 $w\ii v \rvr v_0 u_1\ii $, and $u_0\ii v_0 \rvr v_2 
 u_2\ii $;} \hfill\\
 \text{then there exist $u''$, $v''$, $w''$
 in~$\SS^*$ satisfying}\hfill\\
 \qquad \text{$u\ii v \rvr v''
 {u''}\ii $, $v_1\ii v'' \rvr v_2 {w''}\ii $, and ${u''}\ii
 u_1 \rvr w'' u_2\ii $.}\hfill
 \end{matrix}
 \end{equation}
 (Conditions~\eqref{E:rcuv} and
 \eqref{E:rcuw} might make sense in a
 non-complemented case would the current
 relation~$\rvr$ be replaced with the extended
 relation~$\rvrq$ of Remark~\ref{R:word}.)
\end{rema}

\begin{rema}
 If $(\SS, \RR)$ is an $r$-complete complemented
 presentation, then $r$-rever\-sing is
 compatible with~$\eqp$ in the sense that, if we have
 $u\ii v \rvr v' {u'}\ii $ and $u_1 \eqp u$, then we
 have ${u'}\ii v \rvr v'_1 
 {u'_1}\ii $ for some words~$u'_1$, $v'_1$ satisfying $u'_1
 \eqp u_1$ and $v'_1 \eqp v_1$. We have no such general
 result here. Indeed, with the previous hypotheses, 
 $r$-completeness gives words $u'_1$,
 $v'_1$, and $w$ satisfying ${u'}\ii v \rvr v'_1 
 {u'_1}\ii $, $u_1 \eqp u'_1 w$ and $v_1 \eqp v'_1 w$, but
 there is no general reason for $w$ to be empty. 
 Let us say that two words $u_0$, $v_0$ are co-prime
 if the conjunction of $u_0 \eqp u'_0 w_0$ and $v_0
 \eqp v'_0 w_0$ implies $w_0 = \e$. Then, we could
 deduce $w = \e$ above if we knew that $u_1$ and
 $v_1$ are co-prime, \ie, that reversing
 always produces co-prime words. This is
 true in the complemented case, but, not in general, 
 even if $u$ and $v$ are co-prime for each relation
 $\ss u = \tt v$ in~$\RR$, as shows the example
 developed in Remark~\ref{R:comu} below.
\end{rema}

\section{Recognizing completeness}\label{S:crit}

The characterizations of completeness given in
Prop.~\ref{P:cocu} all are infinitary, in that they
involve checking some condition on infinitely many
words. They therefore give us no effective criterion for
proving completeness. We shall establish now such a
criterion in the case of certain presentations called
homogeneous and connected with Noetherianity.

\begin{defi}\label{D:homo}
 We say that a positive presentation~$(\SS, \RR)$ is
 {\it $r$-homogeneous} if the associated
 congruence~$\eqp$ preserves some $r$-pseudolength,
 the latter being defined as a map~$\ll$ of~$\SS^*$ to
 the ordinals satisfying, for every~$\ss$ in~$\SS$ and every~$u$
 in~$\SS^*$,
 \begin{equation}\label{E:pse2}
 \ll(\ss u) > \ll(u).
 \end{equation}
 We say that $(\SS, \RR)$ is {\it homogeneous} if it 
 preserves both an $r$-pseudolength and an
$l$-pseudolength, the latter defined by the
symmetric condition $\ll(u \ss) > \ll(u)$.
\end{defi}

By definition, the congruence~$\eqp$ associated with
a presentation~$(\SS, \RR)$ is the equivalence relation
generated by the pairs $(uvw, uv'w)$ such that $v =
v'$ is a relation of~$\RR$, so saying that $\eqp$ 
preserves~$\ll$ is equivalent to saying that we have
\begin{equation}
 \ll (uvw) = \ll (uv'w)
 \mbox{\quad for~$v = v'$ in~$\RR$
 and~$u$, $w \in \SS^*$} \label{E:pse1}
\end{equation}

If all relations in~$\RR$ consist of words of equal
length, then the length is both an $r$-
and an $l$-pseudolength, and the presentation is
homogeneous. However, completely different types
exist, as the following examples show.

\begin{exam}\label{X:noet}
 The presentation $(a, b \; a b a = b^2)$ is
 homogeneous. Indeed, the mapping~$\ll$ defined by
 $\ll(a) = 1$, $\ll(b) = 2$, and $\ll(uv) = \ll(u) +
 \ll(v)$ is both an $r$- and an $l$-pseudolength.
 
 A slightly more complicated example is $(a, b, c \;
 ab = bac, ac = ca, bc= cb)$, 
 a presentation for the Heisenberg group. Here, no
 function~$\ll$
 satisfying $\ll(uv) = \ll(u) + \ll(v)$ may be a pseudolength. 
 However, if we define $\ll(u)$ to be the length of~$u$
 augmented by the number of pairs~$(i, j)$ with $i < j$ such
 that the $i$-th letter of~$u$ is~$a$ and the $j$-th letter
 is~$b$---so, for instance, we have $\ll(ab) 
 = \ll(bac) = 3$---then $\ll$ is an $r$- and an
 $l$-pseudolength, and the
 presentation is homogeneous.

 Finally, the presentation $(a, b \; ab = a)$ is
 $r$-homoge\-neous, as shows the
 $r$-pseudo\-length~$\ll$
 defined by $\ll(a) = 1$, $\ll(b) = \o$, and $\ll(uv) =
 \ll(v) + \ll(u)$. As the monoid $\MO(a, b; ab =
 a)$ is not $l$-Noetherian since we have $a \divl
 a$, the next result shows that this presentation is
 not homogeneous.
\end{exam}
 
\begin{prop}\label{P:hono}
 The monoid $\MO(\SS; \RR)$ is $r$-Noetherian (\resp
Noetherian)
 if and only if the presentation $(\SS, \RR)$ is
 $r$-homogeneous (\resp homogeneous).
\end{prop}

\begin{proof}
 If $\ll$ is an $r$-pseudolength on~$\SS^*$, it induces a well
 defined mapping~$\cl\ll$ on $\MO(\SS; \RR)$ such that, by definition, $\xx \divr \yy$ implies
 $\cl\ll(\xx) < \cl\ll(\yy)$. Since the ordinals are well
 ordered, the relation~$\divr$ may have no infinite
 descending chain.
 
 Conversely, assume that $M$ is an $r$-Noetherian
 monoid and $(\SS, \RR)$ is a presentation for~$M$.
 Standard arguments of basic set theory (see for
 instance
 \cite{Lev}) show that there exists a map $\r$ of~$M$
 to the ordinals such that $\xx \divl \yy$ implies $\r(\xx) <
 \r(\yy)$. Then the map~$\ll$ defined by~$\ll(u) =
 \r(\cl u)$ is an $r$-pseudolength
 on~$\SS^*$.
\end{proof}

Our main result now is that, when a presentation
is $r$-homogeneous, then, in order to
prove that the presentation is $r$-complete, it is
sufficient to establish the $r$-cube condition for all
triples of {\it letters}.

\begin{prop}\label{P:loct}
 An $r$-homogeneous positive presentation is
 $r$-complete if and only if any one of the following
 equivalent conditions is satisfied:
 
 (i) The strong $r$-cube condition is satisfied
 on~$\SS^*$;
 
 (ii) The strong $r$-cube condition is satisfied
 on~$\SS$;
 
 (iii) The $r$-cube condition is satisfied on~$\SS$.
\end{prop}

We have already seen in Prop.~\ref{P:cocu} that
$r$-completeness is equivalent to~(i), it is clear that (i)
implies~(ii), and we have observed that the strong
$r$-cube condition always implies the $r$-cube
condition, so (ii) implies~(iii). So, we are left with
the question of proving that (iii) implies say
$r$-completeness, which is the nontrivial point. The
argument will be splitted into several intermediate
statements. Until the end of the proof, we assume that
$(\SS, \RR)$ is an $r$-homogeneous presentation, and
we wish to establish $r$-completeness for every
$4$-tuple of words, \ie, we wish to prove that, if $u v'
\eqp v u'$ holds, then there exist some words $u''$,
$v''$, $w''$ satisfying $u\ii v \rvr v'' {u''}\ii $, $u' \eqp u''
w''$, and $v' \eqp v'' w''$. We fix an $r$-pseudolength~$\ll$
on~$\SS^*$ which is invariant under~$\eqp$.

\begin{lemm}\label{L:ind1}
 The $r$-completeness condition holds for all
 $u$, $v$, $u'$, $v'$ satisfying $\ll(u v')~=~0$.
\end{lemm}

\begin{proof}
 The only possibility is $u = v' = v = u' = \e$, and
 taking $u'' = v'' = w'' = \e$ gives the result.
\end{proof}

\begin{lemm}\label{L:ind2}
 Assume that the $r$-cube condition
 holds on~$\SS$, and $r$-complete\-ness holds for
 all $u$, $v$, $u'$, $v'$ with $\ll(u v') < \alpha$.
Then $r$-completeness holds for all $u$, $v$,
 $u'$, $v'$ with with $u$,
 $v \in \SS$ and $\ll(u v') \le \alpha$.
\end{lemm}

\begin{figure}
 \includegraphics{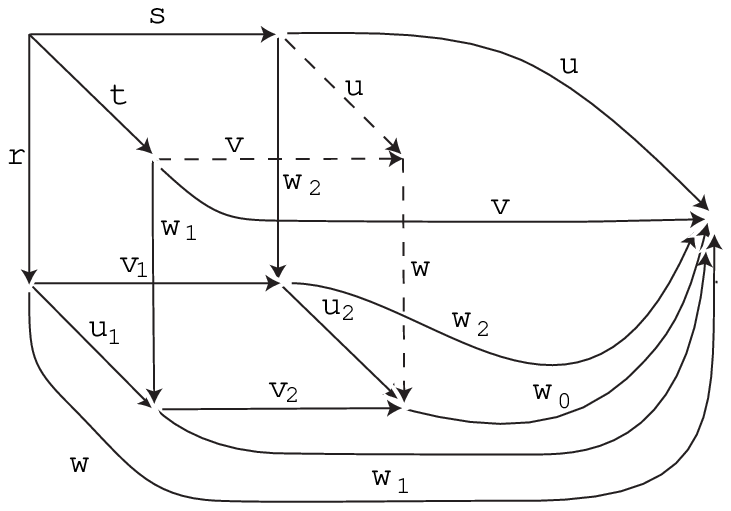}
 \caption{Completeness on~$\SS$}\label{F:indu}
\end{figure}

\begin{proof}
 Assume $\ss v' \eqp \tt u'$ with $\ss$, $\tt \in \SS$ and
 $\ll(\ss u') = \alpha$. We use induction on the minimal
 number of relations~$k$ needed to transform~$\ss v'$
 into~$\tt u'$. The case $k = 0$ corresponds to $\ss v'
 = \tt u'$, hence $\ss = \tt$ and $u' = v'$. In this
 case, taking $u'' = v'' = \e$, $w'' = u'$ gives the
 result. The case $k = 1$ subdivides into two subcases. Either
 the relation connecting $\ss v'$ to~$\tt u'$ does not
 involve the initial letters: then we have $\ss = \tt$, and
 $u' \eqp v'$, and taking $u'' = v'' = \e$, $w''= 
 u'$ gives the result. Or the relation connecting $\ss v'$
 to~$\tt u'$ involves the initial letters: this means that 
 there exists a relation $\ss v'' = \tt u''$ in~$\RR$ and
 a word~$w''$ satisfying $u' = u'' w''$, and $v' =
 v'' w''$: these words $u''$, $v''$, $w''$ give the
 result. 
 
 Assume now $k \ge 2$, and let $\rr w'$ be an
 intermediate word in a shortest path from $\ss v'$ to
 $\tt u'$ (Fig.~\ref{F:indu}). We have $\ss v' \eqp
 \rr w'$ with less
 than $k$~relations, so the induction hypothesis gives
 words $u_1$, $w_1$ and $w'_1$ satisfying
 $\ss\ii \rr \rvr w_1 u_1\ii $, $v' \eqp w_1
 w'_1$, and $w' \eqp u_1 w'_1$. Similarly, we
 have $\rr w' \eqp \tt u'$ with less
 than $k$~relations, so the induction hypothesis gives
 words $v_1$, $w_2$, $w'_2$ satisfying
 $\rr\ii \tt \rvr v_1 w_2\ii $, $w' \eqp v_1
 w'_2$, and $u' \eqp w_2 w'_2$. Then, we have
 $u_1 w'_1 \eqp v_1 w'_2$, and, by definition
 of an $r$-pseudolength, $\ll(u_1 w'_1) < \ll(\rr u_1 w'_1)
 = \ll(\ss w_1 w'_1) = \ll(\ss v') = \alpha$. Applying the
 hypothesis to $u_1$, $v_1$, $w'_1$, $w'_2$, we
 obtain three words $u_2$, $v_2$ and $w'_0$ satisfying
 $u_1\ii v_1 \rvr v_2 u_2\ii $, $w'_1 \eqp v_2 w'_0$, and
 $w'_2 \eqp u_2 w'_0$. At this point, we have $\ss\ii \rr
 \rr\ii \tt \rvr w_1 v_2 u_2\ii w_2\ii $, so the
 hypothesis that $(\SS, \RR)$ satisfies the $r$-cube
condition on
 $\{\ss, \tt , \rr\}$ gives three words
 $u''$, $v''$, $w''_0$ in~$\SS^*$ satisfying $\ss\ii \tt \rvr
 v''{u''}\ii $, $w_1 v_2
 \eqp v'' w''_0$, and $w_2 u_2 \eqp u''
 w''_0$. Put $w'' = w''_0 w'_0$. Then we have
 $u' \eqp w_2 u_2 w'_0 \eqp
 u'' w''$, and $v' \eqp w_1 v_2
 w'_0 \eqp v'' w''$, so the words $u''$,
 $v''$, and $w''$ give the expected result.
\end{proof}

\begin{lemm}\label{L:ind3}
 Assume that $r$-completeness holds for all $u$, $v$,
 $u'$, $v'$ with $\ll(u v') < \alpha$, and for all $u$, 
 $v$, $u'$, $v'$ with $u$, $v \in \SS$ and $\ll(u v')
\le \alpha$. Then $r$-com\-pleteness holds for all $u$, $v$,
 $u'$, $v'$ with $\ll(u v') \le \alpha$.
\end{lemm}

\begin{figure}
 \includegraphics{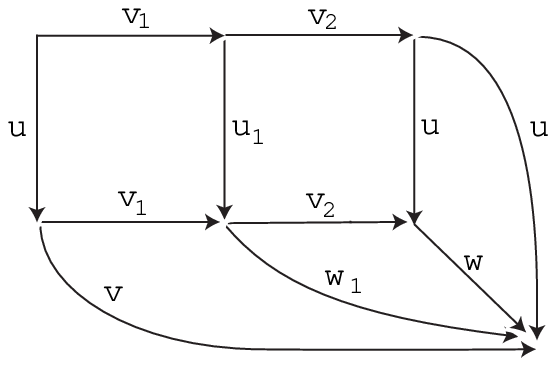}
 \caption{Completeness on~$\SS^*$}\label{F:indv}
\end{figure}

\begin{proof}
 Assume $u v' \eqp v u'$ with $\ll(u v') = \alpha$.
 We wish to prove that there exist $u''$, $v''$, $w''$
 satisfying $u\ii v \rvr v'' {u''}\ii $, $u' \eqp u''
 w''$, and $v' \eqp v'' w''$. If either $u$ or $v$ is empty, the
 result is obvious as, for $u = \e$, we can take $u'' = \e$, $v''
 = v$, and $w'' = u'$. Now, we prove using induction
 on~$m$ that the result holds for $\lg(u) + \lg(v) \le
 m$. By the previous remark, the first nontrivial case is $m =
 2$ with both $u$ and $v$ in~$\SS$. Then the
 conclusion is our second hypothesis. Assume now $m \ge 3$,
 and $v$, say, has length at least~$2$. We write $v = v_1
 v_2$ with both $v_1$ and $v_2$ nonempty
 (Fig.~\ref{F:indv}). The
 hypothesis is $u v' = v_1 (v_2 u')$ with $\ll(u
 v') = \alpha$ and $\lg(u) + \lg(v_1) < m$. Applying the
 induction hypothesis to $u$, $v_1$, $v_2u'$, $v'$,
 we obtain three words
 $u''_1$, $v''_1$, and $w''_1$ in~$\SS^*$ satisfying
 $u\ii v_1 \rvr {u''_1}\ii v_2$, $v_2 u' \eqp u''_1 w''_1$, and
 $v' \eqp v''_1 w''_1$. Now, we have $\ll(v_2 u') < \ll(v_1
 v_2 u') = \alpha$, so applying the first hypothesis to
 $u''_1$, $v_2$, $u'$, $w''_1$, we obtain three words
 $u''$, $v''_2$, and $w''$ satisfying
 ${u''_1}\ii v_2 \rvr v''_2 {u''}\ii $, $u' \eqp u'' w''$, and
 $w''_1 \eqp v''_2 w''$. Put $v'' = v''_1 v''_2$. By
 construction, we have $u\ii v \rvr v''_1 {u''_1} v_2 \rvr v''_1
 v''_2 {u''}\ii $, hence $u\ii v \rvr v'' {u''}\ii $, and
 we have $v' \eqp v''_1 w''_1 \eqp v''_1 v'_2
 w'' = v'' w''$, the expected result.
\end{proof}

It is now easy to complete the proof of 
Prop.~\ref{P:loct}.

\begin{proof}[Proof of Prop.~\ref{P:loct}]
 Assume that $r$-completeness
 fails for some $u_0$, $v_0$, $u_0'$, $v_0'$. Let
 $\alpha$ be the
 minimal possible value of~$\ll(u_0 v_0')$ for such a
 counter-example. By Lemma~\ref{L:ind1},
 $\alpha$ is not~$0$. Now, by construction, the presentation is
 $r$-complete for all $u$, $v$, $u'$, $v'$ with
 $\ll(u v') < \alpha$, hence, by Lemma~\ref{L:ind2}, it
is
 $r$-complete for for all $u$, $v$, $u'$, $v'$ with
 $\ll(u v') \le \alpha$ and $u, v \in \SS$, hence, by
 Lemma~\ref{L:ind3}, it is also $r$-complete for
 all $u$, $v$, $u'$, $v'$ with $\ll(u v') \le \alpha$,
 contradicting the definition of~$\alpha$.
\end{proof}

We can also complete the proof of
Prop.~\ref{P:exno}.

\begin{proof}[End of proof of Prop.\ref{P:exno}]
 We assume that $M$ is a Noetherian monoid, and
 $\RR$ is an $r$-selector on~$\SS$ in~$M$. 
 We have already seen that $(\SS, \RR)$ is a
 presentation for~$M$, and we wish to prove that this
 presentation is $r$-complete. As $M$ is Noetherian, we
 may use a pseudolength~$\ll$, and use the same
 inductive scheme as for Prop.~\ref{P:loct}. Here, by
 definition of an $r$-selector, the
 presentation is $r$-complete at $u$, $v$, $u'$, $v'$
 whenever $u$ and $v$ are single letters, \ie, the
 conclusion of Lemma~\ref{L:ind2} is true directly.
 Then it suffices to use Lemma~\ref{L:ind3} for going
 from $\ll(u v') < \alpha$ to $\ll(u v') \le \alpha$ for
 every~$\alpha$, and deducing $r$-completeness
 for all $u$, $v$, $u'$, $v'$. 
\end{proof}

Returning to the framework of this section, we
deduce from Prop.~\ref{P:loct} the following
(necessary and sufficient) criterion for recognizing
$r$-complete presentations:

\begin{algo}\label{A:strong cube}
 Let $(\SS, \RR)$ be an $r$-homogeneous
 presentation. For each triple of
 letters~$\ss$, $\tt$, $\rr$ in~$\SS$:

 (i) Reverse $\ss\ii \rr\rr\ii \tt$ to all
 possible words of the form~$u v\ii $;
 
 (ii) For each $u v\ii $ so obtained, check 
 $\ss u \eqp \tt v$, or, alternatively,
 $(\ss u)\ii (\tt v) \rvr \e$.
 
 Then $(\SS, \RR)$ is $r$-complete if
 and only if the answer at Step~(ii) is always positive.
\end{algo}

The theoretical interest of the previous result is to
show that $r$-completeness, which is a priori a
$\Sigma^0_1$ (\ie, recursively enumerable, cf.
\cite{Mil}) property, actually is a
$\D^0_1$ (\ie, recursive) property in good cases.

\begin{prop}
 Assume that $(\SS, \RR)$ is a finite homogeneous
 presentation
 such that, for some recursive function~$f$, every
 $r$-reversing sequence from a length~$n$~word
 has length $f(n)$ at most. Then for $(\SS, \RR)$ to be $r$-complete is a recursive property.
\end{prop}

\begin{proof}
 Applying Algorithm~\ref{A:strong cube} involves finitely
 many reversing processes, each of which is assumed
 to have a recursively bounded length. So the whole
 process has a recursively bounded length.
\end{proof}

The main interest of the method presumably lies
in its practical tractability. It can be implemented on a
computer easily, and then be used to test concrete
presentations (when the presentation contains several
relations $\ss u = \tt v$ with the same initial
letters~$\ss$ and $\tt$, $r$-reversing is a
non-deterministic process, and checking the
cube condition by hand quickly becomes impossible).
Observe that, for the computer approach, the strong
cube condition is better suited than the cube
condition, as the only pratical way of proving $u \eqp
v$ is to check that $u\ii v$ reverses to the empty word.

 The completeness criterion of
 Proposition~\ref{P:loct}
 applies in particular in the complemented case. In
 this special case, it had already been proved
 in~\cite{Dff} that the satisfaction of
 Condition~\eqref{E:cohe} for $u$, $v$, $w$
 in~$\SS$, which we have seen is similar to the current
 cube condition, is a sufficient condition for
 $r$-completeness.
 
\begin{exam}
 Let us consider the standard presentation of Artin
 braid groups, or, more generally, of any Artin group
 with finite Coxeter type. Then the presentation is
 homogeneous, as all relations preserve the length of
 the words. Then the (strong) cube condition can be
 checked systematically. Observe that it suffices to
 consider the various possible types of relations only.
 For instance, in the case of the braid groups, there
 are only two types of relations, namely the
 length~$2$ relations $\s_i \s_j = \s_j \s_i$ and the
 length~$3$ relations $\s_i \s_j \s_i = \s_j \s_i \s_j$,
 and, therefore, it is sufficient to consider one triple
 of generators for each possible triple of relations,
 so checking the cube condition for the three
 triples $(\s_1, \s_2, \s_3)$ for type $3, 3, 2$, $(\s_1,
 \s_2, \s_4)$ for type $3, 2, 2$, and $(\s_1, \s_3,
 \s_5)$ for type $2, 2, 2$ is enough to claim that the
 standard presentation of {\it every} group~$B_n$ is
 complete. The verification is what Garside makes in
 his proof of Prop.~H in~\cite{Gar}. Similarly,
 the standard presentation of every Artin group is
 complete, as shown in~\cite{BrS}.
 
 More recently, a new presentation of the braid
 group~$B_n$ has been proposed by Birman, Ko, and
 Lee in~\cite{BKL}. This presentation is homogeneous
 and complemented, and the cube condition is
 satisfied, as established in~\cite{BKL}. So the
 presentation is complete, as are more generally
 the so-called dual presentations of the Artin groups
 investigated in~\cite{BDM, PicNote}.
\end{exam}

 Let us mention that other criteria have been
 established subsequently, always in the
 complemented case. In particular, it is proved in
 \cite{Dgk} that, if $(\SS, \RR)$ is a complemented
 presentation (homogeneous or not), then
 the satisfaction of Condition~\eqref{E:cohf} for $u$,
 $v$, $w$ in the closure of~$\SS$ under
$r$-reversing is always a sufficient condition
 for $r$-completeness. This criterion does not seem
 to extend to the general case---nor does either
 the one established in~\cite{Dgc}. The problem
 here is that the cube condition for letters does not
 imply the cube condition for words directly, because
 the elementary cubes cannot be stacked so as to give
 the desired cube. Such an approach can work only if
 we resort to the ``superstrong'' cube
 condition~\eqref{E:rcuv} where all faces are
 reversings. 

\section{Completion of a presentation}
\label{S:completion}

The criterion of Prop.~\ref{P:loct} fails when we find a
cube that cannot be completed using reversing. This
means that some equivalence follows from the relations
of the considered presentation, but that it cannot be
proved using reversing. Now there always exists a way
for forcing some relation $u \eqp v$ to be provable by
reversing, namely adding it to the presentation. Of
course, repairing one obstruction to completeness in
this way may in turn introduce new obstructions.
But we shall see now that the completion process
so sketched always comes to an end, thus yielding a
complete presentation.

Let us begin with an example.

\begin{figure}[htb]
 \includegraphics{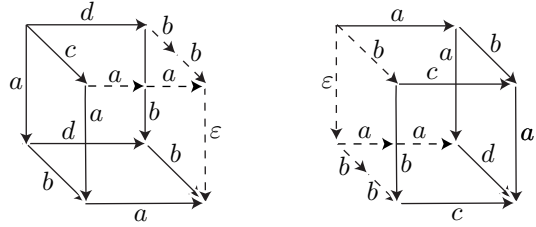}
 \caption{Completion of a presentation}
 \label{F:hako}
\end{figure}

\begin{exam}\label{X:hako}
 (Fig.~\ref{F:hako})
 Let us consider the presentation
 \begin{equation}\label{E:hako}
 (a, b, c, d \; ab=bc = ca, ba = ad = db).
 \end{equation}
 Presentation~\eqref{E:hako} is one of the nonstandard
 presentations of Artin's braid group~$B_3$
 introduced by V.~Sergiescu in~\cite{Ser} and
 considered in~\cite{HaK}: the connection with the
 standard generators~$\s_1$ and~$\s_2$ is
 given by $a = \s_1$, $b = \s_2$, $c =
 \s_1\s_2\s_1\ii$, $d = \s_2\s_1\s_2\ii$. 
 All relations involve words of equal length, so
 \eqref{E:hako} is homogeneous, and Prop.~\ref{P:loct} is
 relevant. Now, when checking the strong cube
 condition for~$(c, a, d)$, we find
 that $c\ii a a\ii d$
 reverses to~$a^2 b^{-2}$, while the presentation
 contains no relation of the form $c \cdots = d
 \cdots$. Here the strong cube condition fails, and
 the presentation~\eqref{E:hako} is not $r$-complete.
 
 The previous failure is due to the relation
 $ca^2 = db^2$, which is a consequence of the
 relations in the presentation, but cannot be 
 proved using reversing associated
 with~\eqref{E:hako}. Now, if we
 add the above relation to the presentation, thus
 obtaining
 \begin{equation}\label{E:hakp}
 (a, b, c, d \; ab=bc = ca, ba = ad = db, ca^2
 = db^2),
 \end{equation}
 then \eqref{E:hakp} is equivalent to ~\eqref{E:hako}
 in that the
 associated monoid and group are the same, and, by
 construction, the relation $ca^2 = db^2$ can now be
 proved by reversing. Of course, new obstructions
 could appear as introducing new relations produces
 new reversing sequences. However, this does not
 happen here, and the reader can
 check that the presentation~\eqref{E:hakp} is $r$-complete.
 
 A symmetric approach is possible for
 $l$-completeness using $l$-reversing and the
 $l$-strong cube condition. The reader can check than the
 presentation~\eqref{E:hakp} is not $l$-complete: we have
 $c a\ii a d\ii \rvl b^{-2} a^2$, and, again, no way
 for proving the relation $a^2 d \eqp b^2 c$ using
 \eqref{E:hakp}-reversing. Once more, the solution is to add
 the missing relation to the presentation, which
 becomes
 \begin{equation}\label{E:hakq}
 (a, b, c, d \; ab=bc = ca, ba = ad = db, ca^2 =
 db^2, a^2 d = b^2 c),
 \end{equation}
 and the reader will now check that \eqref{E:hakq} is
 $l$-complete; it is also $r$-complete as it includes
 \eqref{E:hakp} which is $r$-complete, so, finally,
 \eqref{E:hakq} is a complete presentation.
\end{exam}

The previous example gives a general method for
constructing complete presentations. 

\begin{defi}
 Let $(\SS, \RR)$ be a positive presentation.
 
 \noindent (i) We say that $(\SS, \RR')$ is a {\it
 $1$-completion} of~$(\SS, \RR)$ if there exist
 $\ss$, $\tt$, $\rr$ in~$\SS$ and $u$, $v$
 in~$\SS^*$ satisfying $\RR' = \RR \cup \{\ss v = \tt
 u\}$, $\ss\ii \rr \rr\ii \tt \rvr^{\!R} v u\ii $ but
 $v\ii \ss\ii \tt u \not \rvr^{\!R} \e$.

 \noindent (ii) We say that $(\SS, \RR_\xi)_{\xi <
 \theta}$ is a {\it $r$-completing sequence} if
 $(\SS, \RR_{\xi + 1})$ is a
 $1$-completion of~$(\SS, \RR_\xi)$ for
 each~$\xi$, and, for
 $\xi$ limit, we have $\RR_\xi = \bigcup_{\eta < \xi} \RR_\eta$.
\end{defi}

In other words, the presentation~$(\SS, \RR')$ is a
$1$-completion of~$(\SS, \RR)$ if it is obtained by
fixing one obstruction to the strong $r$-cube condition
for~$(\SS, \RR)$.

\begin{prop}
 Assume that $(\SS, \RR)$ is a homogeneous
 presentation of cardinality~$\kappa$. Then every
 $r$-completing sequence from~$(\SS, \RR)$ ends
 up with an equivalent $r$-complete presentation in
 less than $\sup(\kappa^+, \aleph_1)$~steps.
\end{prop}

\begin{proof}
 Assume first that $(\SS, \RR')$ is a $1$-completion
 for~$(\SS, \RR)$, say $\RR' = \RR \cup \{\ss v = \tt
u\}$.
 By definition, we have $\ss\ii \rr \rr\ii \tt 
 \rvr^{\!R} v u\ii $ for some~$\rvr$, hence, by
 Lemma~\ref{L:weak}(i), $\ss v \eqp^{\RR} \tt u$.
 Therefore, the congruence~$\eqp^{\RR'}$ coincides
 with~$\eqp^{\RR}$, and the presentations
 $(\SS, \RR)$ and $(\SS, \RR')$ are equivalent. 
 Any ($r$-)pseudolength that is preserved
by~$\eqp^{\RR}$ is
 also preserved by~$\eqp^{\RR'}$, so $(\SS, \RR)$
 being ($r$-)homogeneous is equivalent to
 $(\SS, \RR')$ being ($r$-)homogeneous.
 
 If $\SS$ has cardinality~$\kappa$ (finite or
 infinite), then $\SS^*$ has cardinality
 $\sup(\kappa, \aleph_0)$, and so does the set of all
 possible relations
 over~$\SS$. Then the length~$\theta$ of a strictly
 increasing sequence of sets of relations on~$\SS$ say
 $(R_\xi)_{\xi \le \theta}$
 is less than $\sup(\kappa, \aleph_0)^+$, \ie, 
 than $\sup(\kappa^+, \aleph_1)$: otherwise, we
 would obtain an injective mapping~$f$ of the latter
 cardinal into $\SS^* \times \SS^*$ by defining
 $f(\xi)$ to be one element of~$\RR_{\xi +1} \setminus \RR_{\xi}$. The hypothesis that $\RR_\theta$ cannot
 be completed implies that it is $r$-complete.
\end{proof}

In particular, any finite presentation can be
completed in a countable number of steps---but we do
{\it not} claim that, starting from~$(\SS, \RR_0)$ and
defining $(\SS, \RR_{n+1})$ to be a $1$-completion
of~$\RR_n$ implies that $(\SS,
\bigcup_n \RR_n)$ is $r$-complete: the iteration may be
longer than~$\o$. Actually, for practical examples, 
the interesting situation is when the possible
completion requires a finite number of steps only, as
was the case for the presentation of
Example~\ref{X:hako}.

\begin{exam}\label{X:heis}
 Let us consider the standard presentation of the
 Heisenberg group
 \begin{equation}\label{E:heis}
 (a, b, c \;
 ab = bac, ac = ca, bc= cb).
 \end{equation}
 We have seen in Example~\ref{X:noet} that it is
 homogeneous, and, therefore, eligible for our current
 appoach. Now, we find $c\ii b b\ii 
 a \rvr b a b\ii $, but $c\ii a$ only
 reverses to~$a c\ii $, and $b a \eqp a w$, 
 $b \eqp c w$ holds for no word~$w$ on~$\{a,
 b, c\}$.
 According to the scheme above, we add the missing
 relation $c b a =a c$, getting the new presentation
 \begin{equation}\label{E:heit}
 (a, b, c \;
 ab = bac, ac = ca, bc= cb, c b a =a b).
 \end{equation}
 The reader can check that, now, the strong $r$-cube
 condition holds on~$\{a, b, c\}$, and, therefore,
 \eqref{E:heit} is $r$-complete. The latter being
 symmetric, it is actually complete.
\end{exam}

\begin{exam}
 Let us consider the presentation
 \begin{equation}\label{E:bkls}
 (a, b, c \; a^2 = b^2, ab=bc=ca).
 \end{equation}
 One recognizes the Birman-Ko-Lee presentation of the
 braid group~$B_3$, completed with the
 relation~$a^2 = b^2$. Thus, the group
 defined by~\eqref{E:bkls} is the quotient of~$B_3$
 under the relation $\s_1^2 = \s_2^2$. The reader
 can check that \eqref{E:bkls} is not
 $r$-complete, and that completing it leads (in 5
 steps) to the presentation~$(S_3 \; \RR_3)$ of
 Example~\ref{X:aaii}.
\end{exam}

\begin{rema}
 Assume that $u = w$ and $w = v$ are two relations
 in the considered presentation. Then adding the
 relation $u = v$ is a special case of the completion
 procedure described above---which may suggest to
 call {\it transitive} a presentation satisfying the
 cube condition. Indeed, let us isolate the first letters
 in~$u$, $v$, $w$, say $u = \ss u'$, $v = \tt
 v'$ and $w = \rr w'$. Then we have
 \[
 \ss\ii \rr \rr\ii \tt \rvr u' {w'}\ii w' {v'}\ii
 \rvr u' {v'}\ii,
 \]
 and the completion procedure consists in adding the
 relation $\ss u' = \tt v'$, \ie, $u = v$, if we
 cannot obtain $u\ii v \rvr \e$ using the current
 relations. (In the case of Example~\ref{X:heis}, 
 the presentation~\eqref{E:heit} is $r$-complete
 although it contains $ab = bac$ and $ab = cba$
 but not $bac = cba$ because the relation $(bac)\ii
 (cba) \rvr \e$ is already true, and there is no need
 to add $bac = cba$. 
\end{rema}

\section{Reading properties of the monoid}
\label{S:mono}

We enter now the second part of our study. Our aim is to
show that, if $(\SS, \RR)$ is a complete presentation, then 
several properties of the monoid~$\MO(\SS; \RR)$ and
of the group~$\GR(\SS; \RR)$ can be read on the
presentation. We begin with the monoid. We recall
that, when $u$ is a word in~$\SS^*$, then the element of~$\MO(\SS; \RR)$  represented
by~$u$, \ie, the $\eqp$-equivalence class of~$u$, is
denoted by~$\cl u$. 

Let us begin with cancellativity. As mentioned in the
introduction, it is easy to recognize whether a monoid
given by a complete presentation admits cancellation.

\begin{prop}\label{P:canc}
 Assume that $(\SS, \RR)$ is an $r$-complete
 presentation. Then
 the monoid~$\MO(\SS; \RR)$ admits left cancellation if and
 only if $u\ii v \rvr \e$ holds for every relation of the form
 $\ss u = \ss v$ in~$\RR$. In particular, a sufficient condition
 for $\MO(\SS; \RR)$ to admit left cancellation is:
 \begin{equation}\label{E:canr}
 \tag{$C_r$}
 \text{$\RR$ contains no relation $\ss u = \ss v$
 with $u \not= v$.}
 \end{equation}
\end{prop}

\begin{proof}
 The condition is necessary, for $\ss u = \ss v$ belonging
 to~$\RR$ implies $\ss u \eqp \ss v$, hence $u \eqp v$
 if left cancellation is allowed, and, applying
 Prop.~\ref{P:cocu}, $u\ii v \rvr \e$ since
 the presentation is $r$-complete.
 
 Conversely, assume $\ss u \eqp \ss v$ with $\ss \in \SS$. By
 Prop.~\ref{P:cocu}, $u\ii \ss\ii \ss
 v \rvr \e$ holds. By Lemma~\ref{L:deco}, this means that
 there exist words $u'$, $u''$, $v'$, $v''$ satisfying 
 \[
 \ss\ii \ss \rvr v' {u'}\ii, 
 u\ii v' \rvr {u''}\ii , 
 {u'}\ii v \rvr v'', 
 \mbox{ and } {u''}\ii v'' \rvr \e. 
 \]
 By Prop.~\ref{P:eqpo}, this
 implies $u \eqp v' u'' \eqp v' v''$ and $v \eqp u' v''$. Thus, if
 $u' \eqp v'$ or, equivalently, ${u'}\ii v' \rvr \e$,
 holds, we deduce $u \eqp v$, \ie, left cancellation is allowed
 in~$\MO(\SS; \RR)$.
\end{proof}

\begin{coro}
 Assume that $(\SS, \RR)$ is a complete presentation.
 Then a sufficient condition for $\MO(\SS; \RR)$ to be
 cancellative is 
 \begin{equation}\label{E:canc}
 \tag{$C$}
 \text{$\RR$ contains no relation $\ss u = \ss v$ or $u \ss =
 v \ss$ with $u \not= v$.}
 \end{equation}
\end{coro}

\begin{exam}
 All presentations we have considered so far satisfy
 Condition~$(C)$, hence the corresponding monoids
 are cancellative. In particular, so is the
 monoid $\MS$ of Example~\ref{X:hako}.
\end{exam}

Let us consider now the word problem for the 
presentation~$(\SS, \RR)$, \ie, the question of
deciding whether two words~$u$, $v$ in~$\SS^*$ represent
the same element of the monoid~$\MO(\SS; \RR)$, \ie, whether
$u \eqp v$ holds. By Prop.~\ref{P:cocu},
if $(\SS, \RR)$ is an $r$-complete presentation, then $u
\eqp v$ is equivalent to $u\ii v
\rvr \e$, \ie, word equivalence is always detected by
$r$-reversing. As was observed in Sec.~\ref{S:redr}, this
need not give a solution for the word problem if we
have no bound on the length of the reversing
sequences. However, Prop.~\ref{P:dewp} gives the
following sufficient condition:

\begin{prop}\label{P:wdpb}
 Assume that $(\SS, \RR)$ is a finite $r$-complete
presentation satisfying
 \begin{equation}\label{E:fini}
 \tag{$F_r$}
 \text{The closure of~$\SS$
 under $r$-reversing is finite.}
 \end{equation}
 Then the monoid~$\MO(\SS; \RR)$ satisfies a
 quadratic isoperimetric inequality, \ie, every
 relation $u \eqp v$ can be established using
 $O((\lg(u) + \lg(v))^2)$ relations of~$\RR$ at most,
 and its word problem is solvable in quadratic time.
\end{prop}

\begin{proof}
 Let $k$ be the supremum of the number of $r$-reversing
 steps needed to reverse $u\ii v$ into $v' {u'}\ii $ for
 $u$, $v$, $u'$, $v'$ in the closure~$\SSr$
of~$\SS$ under $r$-reversing. Prop.~\ref{P:dewp}
 implies that, if $u$, $v$ are words of length~$p$
 and $q$ respectively and $u \eqp v$ holds, then
 $u\ii v$ reverses to~$\e$ in $k pq$~steps at most,
 hence in $O((p + q)^2)$ reversing steps.
 As each reversing step involves at most one relation
 of~$\RR$ (reversing $\ss\ii \ss$ to~$\e$ requires
 none), we conclude that $u \eqp v$ can be proved
 using at most $O((p + q)^2)$ relations of~$\RR$.
\end{proof}

\begin{exam}
 We already observed that Condition~$(F_r)$ applies
 to the monoids of Example~\ref{X:aaii}: the latter
 therefore satisfy a quadratic isoperimetric inequality.
\end{exam}

Let us consider now common (right) multiples.
By Proposition~\ref{P:eqpo}, $r$-reversing
computes common $r$-multiples in the considered
monoid: $u\ii v \rvr v' {u'}\ii$ implies $u v'
\eqp v u'$, so the element of the monoid represented
by~$uv'$ and $vu'$ is a common right multiple
of~$\cl u$ and~$\cl v$. We can therefore expect 
properties involving common $r$-multiples to be easily
recognized using $r$-reversing.

\begin{prop}
 Assume that $(\SS, \RR)$ is an $r$-complete
presentation. Then a necessary and sufficient 
 condition for any two elements of~$\MO(\SS; \RR)$ to
 admit a common right multiple is
 \begin{equation}\label{E:exis}
 \tag{$E_r$}
  \begin{matrix}
   \text{There exists $\SSp$ satisfying
   $\SS \ince \SSp \ince \SS^*$ and} \hfill\\
   \quad \text{for all $u$, $v$ in~$\SSp$, 
   there exist $u'$, $v'$ in~$\SSp$ 
   satisfying $(u v')\ii (v u') \rvr \e$.}\hfill
  \end{matrix}
 \end{equation}
\end{prop}

\begin{proof}
 Assume that any two elements of~$\MO(\SS; \RR)$
 admit a common right multiple. This means that,
 for all words~$u$, $v$ in~$\SS^*$, there exist
 two words~$u'$, $v'$ satisfying $u v' \eqp v u'$,
 \ie, equivalently, $(u v')\ii (v u') \rvr \e$, since
 $(\SS, \RR)$ is $r$-complete. So $\SSp = \SS^*$
 is convenient.
 
 Conversely, assume that $\SSp$ satisfies
 Condition~$(E_r)$. The latter implies that, for all $u$,
 $v$ in~$\SS'$, there exist $u'$, $v'$ in~$\SS'$
 satisfying $u v' \eqp v u'$. Then, an easy
 induction on~$p + q$ shows that, for $u$
 in~${\SSp}^p$ and $v$ in~${\SSp}^q$,
 there exist $u'$ in~${\SSp}^p$ and $v'$ in~${\SSp}^q$
 satisfying $uv' \eqp vu'$, and, therefore, the elements
  of~$\MO(\SS; \RR)$ represented by~$u$ and~$v$
  admit a common $r$-multiple.
\end{proof}

\begin{figure}
 \includegraphics{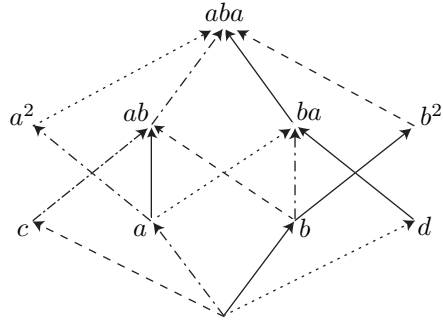}
 \caption{Cayley graph of Sergiescu's monoid~$\MS$}
 \label{F:cayl}
\end{figure}

\begin{exam}
 Let us consider again the monoid~$\MS$ of
 Example~\ref{X:hako}. As shown in Fig.~\ref{F:cayl},
 the family $\{1, a, b, c, d, a^2, ab, ba, b^2, aba\}$
 has the desired properties. So the monoid~$\MS$ 
 admits common right multiples.
\end{exam}

\begin{rema}
 We may replace the relation $(u v')\ii (v u') \rvr \e$
 in Condition~$(E_r)$ by $u\ii v \rvr v' {u'}\ii$, but
 the resulting condition~$(E'_r)$ is stronger, and,
 therefore, more difficult to check in practice. Indeed, 
 $u\ii v \rvr v' {u'}\ii$ implies $u v' \eqp v u'$, hence
 $(u v')\ii (v u') \rvr \e$ for an $r$-complete
 presentation, so $(E'_r)$ implies~$(E_r)$. But,
conversely, $(u
  v')\ii (v u') \rvr \e$ implies that $u\ii v \rvr
 v''{u''}\ii$ holds for some words~$u''$,~$v''$, but the
 hypothesis that $u'$, $v'$ can be chosen in~$\SS'$
 need not imply that $u''$, $v''$ do. 
\end{rema}

As for the existence of least common multiples, we have the
following criterion:

\begin{prop}\label{P:lcmu}
 Assume that $(\SS, \RR)$ is an $r$-complete
presentation. Then a sufficient condition for any two
 elements of~$\MO(\SS; \RR)$ admitting a common
 right multiple to admit a least one is that $(\SS, \RR)$
 is
 an $r$-complemented presentation, \ie, it satisfies
 the condition
 \begin{equation}\label{E:uniq}
 \tag{$U_r$}
   \begin{matrix}
     \text{$\RR$ contains no relation $\ss u = \ss v$
      with $u \not= v$, and,}\hfill\\
      \quad\text{for $\ss \not= \tt$,
      it contains at most one relation
      $\ss \cdots = \tt \cdots$.}\hfill
    \end{matrix}
 \end{equation}
 In this case, the $r$-lcm of~$\cl u$ and~ $\cl v$
 is~$\cl{uv'}$, where $u'$ and $v'$ are the unique
 words satisfying $u\ii v \rvr v' {u'}\ii $.
\end{prop}

\begin{proof}
 If the presentation~$(\SS, \RR)$ is complemented,
 $r$-reversing is a
 deterministic process, so, for every pair of words~$u$, $v$
 in~$\SS^*$, there exists at most one pair of
 words~$u''$, $v''$ in~$\SS^*$ satisfying $u\ii v \rvr v''
 {u''}\ii $. Assume that $uv'$ and $vu'$ represent some
 common right multiple of~$\cl u$ and $\cl v$ in~$\MO(\SS; \RR)$. Then, by definition of $r$-completeness, there must
 exist~$w$ satisfying $u' \eqp u'' w$ and $v' \eqp v''
 w$, where $(u'', v'')$ is the unique pair satisfying
 $u\ii v \rvr v'' {u''}\ii $: this means that $\cl{u
 v''}$ is a right lcm of~$\cl u$ and~$\cl v$.
\end{proof}

\begin{exam}
 The criterion applies to the standard or
 dual presentations of the (generalized) braid groups,
 and to the many examples of~\cite{Pic}, so, in
 each case, elements of the associated monoids that
 admit common multiples admit lcm's---as was already
 observed in previous papers dealing with reversing
 in the complemented case. In contradistinction, none
 of the presentations considered in
 Sec.~\ref{S:completion} is complemented, and it is
 easy to check that lcm's do not exist there.
 
 Observe that $r$-completeness is needed for 
 Condition~$(U_r)$ to imply anything. For
 instance, $(U_r)$ is true for the
 presentation~\eqref{E:heis} of the Heisenberg 
 monoid of Example~\ref{X:heis}, though $a$ and
 $c$ have no $r$-lcm in the Heisenberg monoid:
 indeed, $ac$ and $bac$ are distinct $r$-mcm's of~$a$
 and~$c$, but neither is a multiple of
 the other. Now, of course, $(U_r)$ fails for the
 $r$-complete presentation~\eqref{E:heit}.
\end{exam}

\begin{rema}\label{R:comu}
 Prop.~\ref{P:lcmu} tells us that, in the
 complemented case, $r$-reversing computes
 $r$-lcm's, and we could expect that, in the general
 case, it computes $r$-mcm's (minimal common
 multiples). This need not be the case, even for a
 homogeneous presentation. It is true that, if
 $(\SS, \RR)$ is an $r$-complete presentation, then
 every possible $r$-mcm of~$\cl u$
 and $\cl v$ in~$\MO(\SS; \RR)$ can be represented
 by~$uv'$ and~$vu'$ such that  $u\ii
 v \rvr v' {u'}\ii $ holds. Indeed, if $\cl{uv'}$ is an
$r$-mcm
 of~$\cl u$ and~$\cl v$, then $r$-completeness gives
 $u''$, $v''$, $w$ satisfying $u\ii v \rvr v'' {u''}\ii $, $u'
 \eqp u'' w$, and $v' = v'' w$, and the minimality
 of~$\cl{uv'}$ implies that $w$ must be empty.
 But, conversely, it is not true in general that $u\ii
 v \rvr v' {u'}\ii$ implies that $uv'$ and $vu'$ represent
 an $r$-mcm of~$\cl u$ and~$\cl v$, as shows the
 following counter-example: We have seen that
  the presentation $(a, b \; ab=ba, a^2 =
  b^2)$ is homogeneous and complete. Moreover, each
  relation represents an $r$-mcm. However, we have
  $a\ii b^2 \rvr b^2a\ii $, but $ab^2$ is not an $r$-mcm of~$a$ and~$b^2$ as
  $a$ is a common right divisor of $b^2$ and $a$.
\end{rema}

\section{Reading properties of the group}
\label{S:embe}

We turn to the question of reading properties of the
group~$\GR(\SS; \RR)$ when $(\SS, \RR)$ is a
complete positive presentation. Here we shall
consider the question of whether  $\GR(\SS; \RR)$ is a
group of fractions, and, in this case, study its word
problem.

Recognizing whether $\GR(\SS; \RR)$ is a group of
fractions of the monoid~$\MO(\SS; \RR)$ is easy.
Indeed, it is well-known \cite{CPr} that this happens if
and only if $\MO(\SS; \RR)$ satisfies Ore's conditions,
\ie, it is cancellative and every two elements admit a
common multiple. By gathering results from
Sec.~\ref{S:mono}, we obtain directly:

\begin{prop}\label{P:embe}
 Assume that $(\SS, \RR)$ is a complete presentation.
 Then sufficient conditions for the monoid~$\MO(\SS; \RR)$
 to embed in a group of fractions are
 \begin{gather*}\label{E:xxxx}
 \tag{$C$} 
 \text{$\RR$ contains no relation $\ss u = \ss v$
 or $u \ss = v \ss$ with $u \not= v$,}\\
 \tag{$E_r$}
 \begin{matrix}
 \text{There exists $\SSp$ satisfying
 $\SS \ince \SSp \ince \SS^*$ and such that,} \hfill\\
 \quad \text{for all $u$, $v$ in~$\SSp$, 
 there exist $u'$, $v'$ in~$\SSp$ 
 satisfying $(u v')\ii (v u') \rvr \e$.}\hfill
 \end{matrix}
 \end{gather*}
\end{prop}

\begin{exam}\label{X:frac}
 Typical presentations eligible for the previous
 criterion are the standard presentations of the
 spherical Artin groups, \ie, those associated with a
 finite Coxeter group, or, more generally, all
 presentations of Gaussian groups investigated
 in~\cite{Dfx, Dgk, Pic}. All these presentations
 are complemented.
 
 Now, also eligible are the presentations
 considered in Examples~\ref{X:aaii}, \ref{X:hako},
 and~\ref{X:heis}. In each case, 
 the conditions~$(C)$ and~$(E_r)$ are satisfied,
 and the associated monoid embeds in a group of
 fractions. This holds in particular for Sergiescu's
 monoid~$\MS$ of Example~\ref{X:hako}, of which the
 associated group of
 fractions is the braid group~$B_3$: we thus obtain
 a new decomposition of~$B_3$ as a group of
 fractions, besides the classical decomposition
 associated with the monoid~$B_3^+$ and the
 Birman-Ko-Lee decomposition of~\cite{BKL}
 (this answers a question of~\cite{HaK}).
\end{exam}

Under the hypotheses of Prop.~\ref{P:embe},the
congruence~$\eqp$ that defines the
monoid $\MO(\SS; \RR)$ is the restriction of the
congruence~$\eqpm$ that defines the
group~$\GR(\SS; \RR)$, and standard arguments
then imply that $v u\ii \eqpm v' {u'}\ii$ is true
if and only if there exist~$w$ and~ $w'$ satisfying $u w
\eqp u' w'$ and $v w \eqp v' w$. We shall now
reprove and extend this result by establishing a more
precise connection between the
congruences~$\eqpm$, $\eqp$ and the $r$-reversing
relation in the more general case when only $(C_r)$
and~$(E_r)$ are assumed.

\begin{prop}\label{P:char}
 Assume that $(\SS, \RR)$ is an $r$-complete presentation
 satisfying Conditions~$(C_r)$ and~$(E_r)$. 
 
 (i) For all words~$\ww$, $\ww'$ on~$\SS \cup
 \SS\ii$, 
 the relation $\ww \eqpm \ww'$ is true if and only if
 there exist~$u$, $v$, $w$, $u'$, $v'$, $w'$
 in~$\SS^*$ satisfying
 \begin{equation}\label{E:cong}
 \ww \rvr v u\ii , \quad
 \ww' \rvr v' {u'}\ii , \quad
 u w \eqp u' w', \quad
 v w \eqp v' w'.
 \end{equation}
 
 (ii) In particular, for all words $u$, $u'$ in~$\SS^*$,
 the relation $u \eqpm u'$ is true if and only if there
 exists $w$ in~$\SS^*$ sastifying $u w \eqp u' w$.
\end{prop}

The proof will be splitted into several steps. We
assume until the end of the proof of Prop.~\ref{P:char}
that $(\SS, \RR)$ is an $r$-complete presentation
satisfying~$(C_r)$, and~$(E_r)$. For $\ww$, $\ww'$
words on~$\SS \cup \SS\ii $, we say that $\ww \cong
\ww'$ is true if there exist $u$, $v$, $w$, $u'$, $v'$,
$w'$ satisfying~\eqref{E:cong}. Our aim is to prove
that the relations~$\eqpm$ and~$\cong$ coincide.

\begin{lemm}\label{L:cong}
 Assume $\ww \rvr v u\ii $ and $\ww \rvr
 v' {u'}\ii $ with $u$, $v$, $u'$, $v' \in \SS^*$.
 Then we have $v u\ii \cong v' {u'}\ii $.
\end{lemm}

\begin{figure}[htb]
 \includegraphics{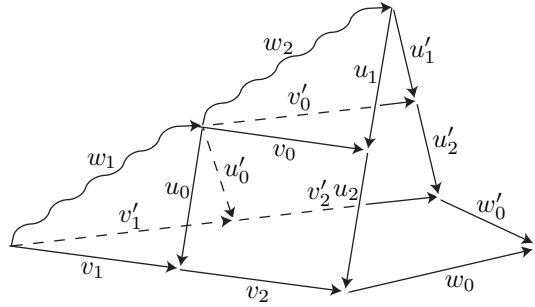}
 \caption{Several reversings}\label{F:cong}
\end{figure}

\begin{proof}
 It suffices to show that there exist two words
 $w$, $w'$ on~$\SS$ satisfying $u w \eqp u'
 w'$ and $v w \eqp v' w'$. The
 hypothesis that $(\SS, \RR)$ satisfies
 $(E_r)$
 implies that there exist words~$w$, $w'$ satisfying
 $v w \eqp v' w'$, and we are left with the 
 question of proving that $v w \eqp v' w'$ implies $u
 w \eqp u' w'$ whenever some word~$\ww$ reverses
 both to $v u\ii $ and to $v' {u'}\ii $. We establish
 the latter implication using induction on the length
 of~$\ww$. The result is trivial if $\ww$ is empty. Assume
 that $\ww$ has length~$1$. If $\ww$ is a letter in~$\SS$,
 say~$\ss$, the hypothesis is $w \eqp w'$, and the
 expected conclusion is $\ss w \eqp \ss w'$, so the
 implication is always true. If $\ww$ is a letter in~$\SS\ii $, 
 the hypothesis is $\ss w \eqp \ss w'$,
 and the expected conclusion is $w \eqp w'$: so
 the implication is true provided $\MO(\SS; \RR)$ admits
 left cancellation. 
 
 Assume now $\ww = \ww_1 \ww_2$ with $\lg(\ww_i) <
 \lg(\ww)$. By Lemma~\ref{L:iter}, there exist
 words $u_i$, $v_i$, $u'_i$, $v'_i$, $i = 0, 1, 2$
 satisfying $\ww_1 \rvr v_1 u_0\ii $, $\ww_2 \rvr v_0
 u_1\ii $ and $u_0\ii v_0 \rvr v_2 u_2\ii $, and
 similar dashed relations (see Fig.~\ref{F:cong}). By
 hypothesis, we have $v_1 v_2 w \eqp v'_1
 v'_2 w'$ and $\ww_1$ reverses both to $v_1
 u_0\ii $ and $v'_1 {u'_0}\ii $, so applying the
 induction hypothesis to~$\ww_1$ gives $u_0 v_2
 w \eqp u'_0 v'_2 w'$, hence $v_0 u_2
 w \eqp v'_0 u'_2 w'$. Now $\ww_2$
 reverses both to $v_0 u_1\ii $ and $v'_0
 {u'_1}\ii $, so applying the induction hypothesis
 to~$\ww_2$ gives $u_1 u_2 w \eqp u'_1
 u'_2 w'$, \ie, $u w \eqp u' w'$, as
 was expected.
\end{proof}

\begin{lemm}\label{L:mult}
 For~$\ss$ in~$\SS$, $\ww \cong \ww'$
 implies $\ww \ss \cong \ww' \ss$ and $\ww \ss\ii \cong
 \ww' \ss\ii $.
\end{lemm}

\begin{figure}[htb]
 \includegraphics{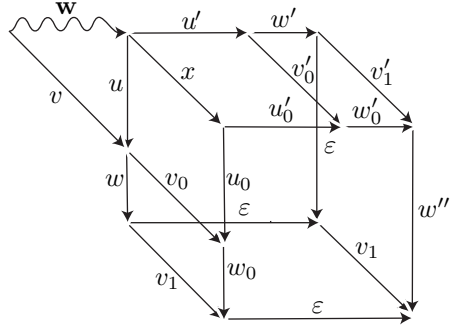}
 \caption{Compatibility with multiplication}
 \label{F:mult}
\end{figure}

\begin{proof}
 Assume
 \[
 \ww \rvr v u\ii , \quad
 \ww' \rvr v' {u'}\ii , \quad
 u w \eqp u' w', \quad
 v w \eqp v' w'.
 \]
 As Condition~$(E_r)$ is satisfied, there exist
 $u_0$, $v_0$, and $v_1$, $w_0$ in~$\SS^*$
 satisfying $u\ii 
 \ss \rvr v_0 u_0\ii $ and $w\ii v_0 \rvr v_1 w_0\ii $
 (Fig.~\ref{F:mult}). So, by construction,
 we have $\ww \ss \rvr (v v_0) u_0\ii $ and $\ss\ii (uw)
 \rvr u_0w_0 v_1\ii $. As the presentation is $r$-complete,
 $u w \eqp u' w'$ implies $(uw)\ii (u'w') \rvr \e$, and, by
 definition, we have $v_1\ii \e \rvr \e v_1\ii $, hence
 $\ss\ii (uw) (uw)\ii (u'w') \rvr u_0 w_0 v_1\ii $.
 The cube condition for $\ss$, $uw$, and $u'w'$
 holds, so there must exist words~$u''$, $v''$, $w''$
 in~$\SS^*$ satisfying 
 $\ss\ii u' w' \rvr u'' {v''}\ii $, $u'' w'' \eqp u_0 w_0$, and $v''
 w'' \eqp v_1$. By
 Lemma~\ref{L:iter}, there exist $u'_0$, $v'_0$, $w'_0$ and
 $v'_1$ satisfying $\ss\ii u' \rvr u'_0 {v'_0}\ii $,
 ${v'_0}\ii w' \rvr w' {v'_1}\ii $, 
 $u'' = u'_0 w'_0$, and $v'' = v'_1$, hence $u_0 w_0 \eqp
 u'_0 w'_0 w''$ and $v_1 = v'_1 w''$. So, we have
 \begin{equation}\label{E:redr}
 \ww\ss \rvr (v v_0) u_0\ii , \quad
 \ww'\ss \rvr (v' v'_0) {u'_0}\ii .
 \end{equation}
 Now we check
 \begin{gather}
 \ss u_0 w_0 \eqp u w v_1 \eqp u' w' v'_1 w'' \eqp 
 \ss u'_0 w'_0 w'', \label{E:10}\\
 v v_0 w_0 \eqp v w v_1 \eqp v' w' v'_1 w''
 \eqp v' v'_0 w'_0 w''.\label{E:11}
 \end{gather}
 As left cancellation is possible, \eqref{E:10} implies
 $u_0 w_0 \eqp u'_0 (w'_0 w'')$, while \eqref{E:11} reads
 $(vv_0) w_0 \eqp (v' v'_0)(w'_0 w'')$, which, together
 with~\eqref{E:redr}, gives $\ww\ss \cong \ww'\ss$.
 
 The case of~$\ss\ii $ is trivial: with the same notation, we
 have
 \[
 \ww\ss\ii \rvr v (\ss u)\ii , \quad
 \ww'\ss\ii \rvr v' (\ss u')\ii , \quad 
 (\ss u) w \eqp (\ss u') w', \quad
 v w \eqp v' w',
 \]
 so $\ww\ss\ii \cong \ww'\ss\ii $ holds as well.
\end{proof}

\begin{proof}[Proof of Prop.~\ref{P:char}]
 (i) Assume $\ww \cong \ww'$. With the notation
 of~\eqref{E:cong}, we find
 \[
 \ww \eqpm v u\ii \eqpm v w w\ii u\ii \eqpm v' w' {w'}\ii {u'}\ii 
 \eqpm v' {u'}\ii \eqpm \ww',
 \]
 so $\ww \cong \ww'$ implies $\ww \eqpm \ww'$

 Conversely, we shall prove that $\cong$ is a congruence
 that contains pairs generating~$\eqpm$. By definition, the
 relation~$\cong$ is reflexive and 
 symmetric. Assume $\ww \cong \ww' \cong \ww''$. This
 means that there exist words $u$, \pp, $w''$ in~$\SS^*$
 satisfying
 \begin{gather*}
 \ww \rvr v u\ii , \quad
 \ww' \rvr v_1' {u_1'}\ii , \quad
 u w \eqp u_1' w'_1, \quad
 v w \eqp v_1' w'_1, \\
 \ww' \rvr v'_2 {u'_2}\ii , \quad
 \ww'' \rvr v'' {u''}\ii , \quad
 u'_2 w'_2 \eqp u'' w'', \quad
 v'_2 w'_2 \eqp v'' w''.
 \end{gather*}
 By Lemma~\ref{L:cong}, there exist $w_1$, $w_2$ in~$\SS^*$
 satisfying $u'_1 w_1 \eqp u'_2 w_2$ and
 $v'_1 w_1 \eqp v'_2 w_2$. Now, as common right multiples
 exist in the monoid $\MO(\SS; \RR)$, we can find
 $w_3$, $w'_3$, $w_4$, $w'_4$ in~$\SS^*$ satisfying
 $w_1 w_3 \eqp w'_1 w'_3 \eqp w_2 w_4 \eqp w'_2 w'_4$,
 and we find
 \begin{gather*}
 u w w'_3 \eqp u'_1 w'_1 w'_3 \eqp u'_1 w_1 w_3
 \eqp u'_2 w_2 w_4 \eqp u'_2 w'_2 w'_4
 \eqp u'' w'' w'_4\\
 v w w'_3 \eqp v'_1 w'_1 w'_3 \eqp v'_1 w_1 w_3
 \eqp v'_2 w_2 w_4 \eqp v'_2 w'_2 w'_4
 \eqp v'' w'' w'_4,
 \end{gather*}
 so the words $w w'_3$ and $w'' w'_4$ witness for $\ww
 \cong \ww''$. So $\cong$ is an equivalence relation.
 
 We claim now that $\cong$ is a congruence, \ie, it is
 compatible with multiplication on both sides. It suffices
 to consider the case of right of left multiplication by a
 single positive or negative letter. Lemma~\ref{L:mult}
 gives the result for right multiplication, and we observe that
 $\ww \cong \ww'$ is equivalent to $\ww\ii \cong {\ww'}\ii $,
 so the result for left multiplication follows.
 
 By definition, $\eqpm$ is the congruence on~$(\SS \cup
 \SS\ii)^*$ generated
 by all pairs~$\{u, v\}$ in~$\RR$, completed with all pairs
 $\{\ss\ss\ii , \e\}$ and $\{\ss\ii \ss, \e\}$ with $\ss \in\SS$.
 Writing
 \begin{gather*}
 u \rvr u, \quad v \rvr v, \quad \e \e \eqp \e \e, \quad 
 u \e \eqp v \e,\\
 \ss\ss\ii \rvr \ss\ss\ii , \quad \e \rvr \e \e\ii , 
 \quad \ss\e \eqp \e \ss, \quad \ss\e \eqp \e \ss, \\
 \ss\ii \ss \rvr \e\e\ii , \quad \e \rvr \e \e\ii , 
 \quad \e\e \eqp \e \e, \quad \e\e \eqp \e \e,
 \end{gather*}
 we see that $u \cong v$, $\ss\ss\ii \cong \e$, and $\ss\ii \ss
 \cong \e$ hold, and we conclude that $\eqpm$ is 
 included in~$\cong$, \ie, that $\ww \eqpm \ww'$ implies
 $\ww \cong \ww'$, which completes the proof of~(i).

 (ii) As $\eqp$ is included in~$\eqpm$, the existence of
 a word~$w$ satisfying $u w \eqp u' w$ is a sufficient 
 condition for~$u \eqpm u'$. Conversely, assume $u
\eqpm u'$. By ~(i),
 $u$ and $u'$ have to reverse to fractions
 satisfying~\eqref{E:cong}. As $u$ and $u'$ belong
 to~$\SS^*$, the only possibilities are $u \rvr u \e\ii $ and $u'
 \rvr u' \e\ii $, so \eqref{E:cong} reduces to the existence of
 $w$, $w'$ in~$\SS^*$ that satisfy $u w \eqp u' w'$ and $\e
 w \eqp \e w'$: this implies $u w \eqp u' w$.
\end{proof}

Let us now return to the hypotheses of
Prop.~\ref{P:embe}, \ie, to the case when the
group~$\GR(\SS; \RR)$ is a group of fractions for the
monoid~$\MO(\SS; \RR)$. The following result shows
that the word problem can always be solved by a
double $r$-reversing, or, alternatively, an
$r$-reversing followed with an $l$-reversing.

\begin{prop}\label{P:dbrv}
 Assume that $(\SS, \RR)$ is a complete presentation
 satisfying Conditions~\eqref{E:canc}
 and~$(E_r)$. Then, for every word~$\ww$
 on~$\SS \cup \SS\ii $, the following are equivalent:
 
 (i) We have $\ww \eqpm \e$; 
 
 (ii) There exist $u$, $v$ in~$\SS^*$ satisfying $\ww
\rvr v
 u\ii $ and $u\ii v \rvr \e$; 

 (iii) There exist $u$, $v$ in~$\SS^*$ satisfying $\ww
\rvr v u\ii 
\rvl \e$.
\end{prop}

\begin{proof}
 Under the hypotheses, we know that, for every word~$\ww$,
 there exist positive words~$u$, $v$ satisfying $\ww \rvr v
 u\ii $. Then $\ww \eqpm \e$ is equivalent to $u \eqpm v$, 
 hence to $u \eqp v$ by Prop.~\ref{P:embe}, and,
 therefore, both to $u\ii v \rvr \e$ and to $v u\ii \rvl \e$
 by Prop.~\ref{P:cocu}.
\end{proof}

\begin{prop}\label{P:quis}
 Assume that $(\SS, \RR)$ is a complete presentation
 satisfying Conditions~\eqref{E:fini}, \eqref{E:canc}
 and~$(E_r)$, \ie, 
 \begin{gather*}\label{E:xxxx}
 \tag{$F_r$} 
 \text{The closure of~$\SS$ 
 under $r$-reversing is finite,}\\
 \tag{$C$} 
 \text{The presentation $\RR$ contains no relation 
 $\ss u = \ss v$ or $u \ss = v \ss$ 
 with $u \not= v$,}\\
 \tag{$E_r$}
 \begin{matrix}
 \text{There exists $\SSp$ satisfying
 $\SS \ince \SSp \ince \SS^*$ and
 such that,} \hfill\\
 \quad \text{for all $u$, $v$ in~$\SSp$, 
 there exist $u'$, $v'$ in~$\SSp$ 
 satisfying $(u v')\ii (v u') \rvr \e$.}\hfill
 \end{matrix}
 \end{gather*}
 Then the group~$\GR(\SS; \RR)$ satisfies a quadratic
 isoperimetric inequality.
\end{prop}

\begin{proof}
 We gather Prop.~\ref{P:dbrv}, which reduces the
 word problem in~$\GR(\SS; \RR)$ to a double
 reversing process,
 and Prop.~\ref{P:wdpb}, which gives a bound on the
 complexity of the latter process.
\end{proof}

\begin{exam}\label{X:ispe}
 The previous criterion applies to the groups
 defined by the complemented presentations of
 Example~\ref{X:frac}. But it also applies to the groups
 associated with the presentations of
 Example~\ref{X:aaii}, thus typically to the
 groups
 \begin{gather*}
 (a, b \; a^2 = b^2, ab = ba)\\
 (a, b, c \; a^2 = b^2 = c^2, ab = bc = ca, ac = ba =
 cb).
 \end{gather*}
 (We recall that the latter is the quotient of~$B_3$
 under the additional relation $\s_1^2 = \s_2^2$.)
 These groups therefore satisfy a quadratic
 isoperimetric inequality. So does the group associated
 with the monoid of Example~\ref{X:hako}, but we saw
 that the latter
 group is~$B_3$, and that result is well known.
 
 Let us consider now the Heisenberg
 group~$H$. The closure of $\{a, b, c\}$ under
 $r$-reversing
 with respect to the (incomplete)
 presentation~\eqref{E:heis} is the infinite
 set $\{\e, a, b, c\} \cup \{ac^n \; n \ge 1\}$, and,
 using the latter, we easily conclude that common
 right multiples exist in the associated monoid, of
 which $H$ is the group of fractions. Then
 Prop~\ref{P:dbrv} shows how to solve the word
 problem using a double reversing with respect to
 the complete presentation
 \begin{equation}
 (a, b, c \;
 ab = bac, ac = ca, bc= cb, c b a =a b).
 \end{equation}
 It can be checked that the complexity of the
 procedure is cubic, which could be expected $H$ is
 known to admit a cubic isoperimetric
 function~\cite{Eps}. 
\end{exam}

In the complemented case, the study proceeds farther,
and it is known that, under the hypotheses of
Prop~\ref{P:quis}, the group~$\GR(\SS; \RR)$ is a
Garside group and, in particular, it is
torsion-free\cite{Dfz} and admits a bi-automatic
structure~\cite{Dgk}. The question of whether the
latter result extends to the general case of non
necessarily complemented presentations seems to be
difficult, as the automatic structures known in the
complemented case relie on the uniqueness of the
gcd's. In any case, the answer is connected with the fine
structure of divisibility in the monoid~$\MO(\SS;
\RR)$, and the importance of words and reversing
becomes secondary. So we shall not discuss the
question here, but refer to~\cite{Dgo} where
the question is investigated directly. Let us mention
that the groups of Example~\ref{X:ispe} turn out to be
automatic.

\bigskip

The above study has led to results about the
group~$\GR(\SS; \RR)$ only in the case when the
latter happens to be a group of fractions for the
monoid~$\MO(\SS; \RR)$. The main open question
now is to determine to which extent word reversing can
be used to prove results about the group~$\GR(\SS;
\RR)$ in the general case. In particular, it would be
interesting to know whether reversing techniques can
be used to study the possible embeddability of the
monoid~$\MO(\SS; \RR)$ in the group~$\GR(\SS;
\RR)$. Let us observe here that the presentation
\begin{equation}\label{E:nemb}
 (a, b, c, d, a', b', c', d'; 
 aa' = bb', ca' = db', ac' = bd')
\end{equation}
introduced in~\cite{Kas} is complete and it satisfies
Condition~$(C)$, so the associated monoid is
cancellative, but the latter does not embed in the
corresponding group, as $cc' = dd'$ holds in the group
(we have there $c\ii d = a' {b'}\ii = a\ii b = c'
{d'}\ii$) but not in the monoid (we do not have ${c'}\ii
c\ii d d' \rvr \e$). Can this negative result be read
directly on  Presentation~\ref{E:nemb}? Similarly, but
on the other direction, it is known that every Artin
monoid embeds in the corresponding group \cite{Par},
but the remarkable proof of the result uses an indirect
approach via a linear representation (inspired
by~\cite{Kra}). Could reversing be used here?

We shall conclude this paper with a more precise
question. Assume that $(\SS, \RR)$ is a positive group
presentation, and let $\rv$ denote the union of the
relations~$\rvrq$ and~$\rvlq$, \ie, the extended
$r$-reversing considered in Remark~\ref{R:word} and
its left counterpart. Prop.~\ref{P:dbrv} tells us
that, if $(\SS, \RR)$ is a complete presentation such
that the monoid~$\MO(\SS; \RR)$ is cancellative and
admits common right multiples, then a word~$\ww$
represents~$1$ in the group~$\GR(\SS; \RR)$ if and
only if $\ww \rv \e$ holds. If common multiples do not
exist in~$\MO(\SS; \RR)$, the proof is no longer valid.
However, the above result, namely
that $\ww \eqpm \e$ is equivalent to
$\ww\rvr \e$, also holds in the case of a free group
\ie, when $\RR$ is empty: in this
case, reversing coincides with free reduction, and it is
true that $\ww$ represents~$1$ in a free group if and
only if it freely reduces to~$\e$ (with an unbounded
number of alternations between $r$- and
$l$-reversing, contrary to the case of
Prop.~\ref{P:dbrv} where one alternation is
enough). Similarly, in the case of
Presentation~\eqref{E:nemb}, the key relation $cc' =
dd'$, which we have seen holds in the group but not in
the monoid, can be proved using reversing, \ie,
$(cc')\ii (dd') \rv \e$ holds, as we find
\[
 {c'}\ii c\ii d d'
 \rvr {c'}\ii a' {b'}\ii d'
 \rvl {c'}\ii a\ii b d'
 \rvr {c'}\ii c' {d'}\ii d' \rvr \e
\]
(with two alternations between $r$- and $l$-reversing).
This leads to the general problem of whether the word
problem of the group can be solved using reversing.
Simple counter-examples, such as the presentation
$(a, b, c; ab = ac)$ suggested by S.~Lee, show that some
assumptions have to be satisfied, but the following
question is open:

\begin{ques}
 Let $(\SS, \RR)$ be a complete presentation
 satisfying Condition~$(C)$ (so the monoid~$\MO(\SS;
 \RR)$ is cancellative). Is $\ww \rv \e$ a necessary
 (and sufficient) condition for a word~$\ww$ on~$\SS
\cup \SS\ii$ to represent~$1$ in the group~$\GR(\SS;
 \RR)$?
\end{ques}

A positive answer would imply that we can prove $\ww
\eqpm \e$ by introducing no new factor $\ss \ss\ii $ or
$\ss\ii \ss$, so, in some sense, by always going from
one word to another that is not more complicated
(if not shorter, in general). In this sense, solutions for
the word problem based on word reversing are
reminiscent of Dehn's algorithm for hyperbolic groups,
but their range includes more complicated groups, such
as braid groups, or, more generally, Garside groups
(which admit a quadratic isoperimetric function), or
even more complicated groups like the nilpotent
Heisenberg group (which admits a cubic isoperimetric
function). The underlying question is whether one can
prove that a word~$\ww$ is trivial by remaining not
too far from~$\ww$ in the Cayley graph of the
considered group (a precise meaning was given
in~\cite{Dfo}), and reversing gives a positive answer
for many particular groups. The general case is open,
but we conjecture that completeness is relevant.

\end{document}